\documentclass[11pt]{article}

\usepackage[utf8]{inputenc}

\usepackage{xcolor}

\usepackage{amsfonts}
\usepackage{amssymb}
\usepackage{amsmath}
\usepackage{stmaryrd}
\usepackage{theorem}

\usepackage[pdftex]{hyperref}
\hypersetup{linkcolor=blue,colorlinks=true,citecolor=magenta}

\newcommand{\toto}{\ell}

\newcommand{\dd}{\text{\rm d}}
\newcommand{\boxfill}{\hspace*{\fill}}
\newcommand{\Span}{\operatorname{span}}
\newcommand{\Mat}{\mbox{\rm Mat}}
\newcommand{\TT}{{\cal T}}
\newcommand{\gbar}{\overline{g}}
\newcommand{\FF}{{\mathcal X}^{\perp}}
\newcommand{\GG}{{\cal G}}
\newcommand{\afsuris}{\opna{Aff}(g)/\opna{Isom}(g)}
\newcommand{\indaf}{[\opna{Aff}(g):\opna{Isom}(g)]}

\newcommand{\R}{\mathbb R}
\newcommand{\N}{\mathbb N}
\newcommand{\Q}{\mathbb Q}
\newcommand{\Z}{\mathbb Z}

\newcommand{\T}{\mathbb T}
\newcommand{\hop}{\vskip .2cm\noindent}
\newcommand{\hip}{\vskip .1cm\noindent}

\newcommand{\mrm}{\mathrm}

\newcommand{\oD}{\overline D}
\newcommand{\oU}{\overline U}
\newcommand{\oY}{\overline Y}
\newcommand{\oZ}{\overline Z}
\newcommand{\oR}{\overline R}
\newcommand{\Vbar}{\overline{V}}
\newcommand{\Wbar}{\overline{W}}
\newcommand{\Mbar}{\overline{M}}
\newcommand{\ric}{\operatorname{ric}}
\newcommand{\Ric}{\operatorname{Ric}}

\newcommand{\ombar}{\overline\omega}
\newcommand{\opna}{\operatorname}
\newcommand{\wdt}{\widetilde}
\newcommand{\Id}{\opna{Id}}

\newcommand{\ie}{{\em i.e.\@ }}
\newtheorem{enonce}{}[section]
\newtheorem{fact}[enonce]{Fact}
\newtheorem{de}[enonce]{Definition}
\newtheorem{prop}[enonce]{Proposition}

\newtheorem{add-property}[enonce]{Additional property}
\newtheorem{lem}[enonce]{Lemma}

\newtheorem{te}[enonce]{Theorem}
\newtheorem{cor}[enonce]{Corollary}

\newtheorem{rem-notation}[enonce]{Remark/Notation}
\newtheorem{def-notation}[enonce]{Definition/Notation}
\newtheorem{notation}[enonce]{Notation}
\theorembodyfont{\rmfamily}
\newtheorem{rem}[enonce]{Remark}
\newtheorem{example-rem}[enonce]{Example/Remark}

\newtheorem{consequence-notation}[enonce]{Consequence/Notation}

\newtheorem{recall-rem}[enonce]{Recall/Remark}

\newtheorem{recall/vocabulary}[enonce]{Recall/vocabulary}

\newtheorem{reminder}[enonce]{Reminder}

\newcommand{\XX}{{\cal X}}
\newcommand{\BB}{{\cal B}}

\usepackage{multirow}
\usepackage{hhline}

\newlength{\arraycolsepsauvegardegenerale}
\newlength{\arraycolsepsauvegardetemporaire}
\newlength{\arraycolsepsauvegardetemporairebis}

\newenvironment{nnarray}%
{\setlength{\arraycolsepsauvegardetemporaire}{\arraycolsep}\setlength{\arraycolsep}{0cm}\begin{array}}
{\end{array}\setlength{\arraycolsep}{\arraycolsepsauvegardetemporaire}}

\newenvironment{narray}%
{\setlength{\arraycolsepsauvegardetemporaire}{\arraycolsep}\setlength{\arraycolsep}{0.5\arraycolsepsauvegardegenerale}\begin{array}}
{\end{array}\setlength{\arraycolsep}{\arraycolsepsauvegardetemporaire}}

{\setlength{\arraycolsepsauvegardetemporairebis}{\arraycolsep}
\setlength{\arraycolsep}{\arraycolsepsauvegardegenerale}\begin{array}}
{\end{array}\setlength{\arraycolsep}{\arraycolsepsauvegardetemporairebis}}

\title{Affine transformations and parallel lightlike vector fields on compact Lorentzian $3$-manifolds}
\author{Charles Boubel, Pierre Mounoud}
\date{April 4th.\@ 2014, revised November 24th.\@ 2014}
\begin{document}
\maketitle

\begin{abstract}
We  describe the  compact Lorentzian $3$-manifolds admitting a parallel lightlike vector field. The  classification of  compact Lorentzian $3$-manifolds 
admitting non-isometric affine diffeomorphisms follows, together with the complete description of these morphisms. Such a Lorentzian manifold is in some sense an equivariant deformation of a flat one.
\end{abstract}

\section{Introduction}
Let $(M,g)$ be a compact Lorentzian manifold and $D$ its Levi-Civita connexion. We can associate with it, at least, two groups: its isometry group $\opna{Isom}(g)$ and its affine group \ie the group of transformations preserving $D$, that we denote by $\opna{Aff}(g)$. Of course, $\opna{Isom}(g)$ is a closed subgroup of  
$\opna {Aff} (g)$. 
This article has two purposes:\medskip

-- to understand compact Lorentzian $3$-manifolds $(M,g)$ such that $\afsuris$ is non trivial,\medskip

-- to describe the compact Lorentzian $3$-manifolds admitting a parallel lightlike vector field.\medskip

 These two problems are directly related  because when $\afsuris$ is not trivial then there exists non trivial parallel fields of symmetric  endomorphisms of $TM$. Indeed, if 
$\varphi\in \opna{Aff}(g)\smallsetminus \opna{Isom}(g)$ then $\varphi^*g\neq g$ and there exists a field $E$  of symmetric endomorphisms of $TM$, such that $\varphi^*g(.,.)=g(E.,.)$.
As $\varphi$ preserves $D$, the  Levi-Civita connection of the metric $\varphi^*g$ is also $D$ and $E$ is parallel. In the three dimensional situation it means that, up to a finite cover, $(M,g)$ posseses a non trivial parallel vector field $X$.  The cases where $X$ is lightlike quickly appears to be the only interesting one (see Proposition \ref{decomposable}).

Moreover, both problems are interesting in themselves. When a pseudo-Riemannian manifold $(M,g)$ admits a parallel {\em non degenerate} distribution, it splits into a Riemannian product ---~at least locally, and globally as soon as it is complete and simply connected~---, this is the well-known generalization \cite{wu67} by Wu of a theorem of de Rham. When the distribution is degenerate, the situation becomes much richer, and more complicated. For Lorentzian manifolds, the local situation is now quite well understood, see {\em e.g.\@} the survey \cite{galaev-leistner2008}, whereas the subject remains widely open for general pseudo-Riemannian metrics. Some recent works study the {\em global} behaviour of such manifolds. The case of homogeneous spaces, which is quite specific, has been  a little studied, and in the general case, results and examples on the full holonomy group of Lorentzian metrics are also given in \cite{larz-these, baum-larz-leistner}, some of them involving compact manifolds (see {\em e.g.\@} in \cite{larz-these} metrics on circle bundles similar to those we see appearing in our work). However, very few works treat specifically the compact case ; see {\em e.g.\@} \cite{Leistner-Schliebner13} for the case of pp-waves (compact implies complete), \cite{Schliebner12} with a parallel Weyl tensor or \cite{Schliebner13} with an assumption on the first Betti number. This work is devoted to it, more precisely to compact Lorentzian manifolds with a parallel, lightlike vector field, in the lowest interesting dimension {\em i.e.\@} three, but with no additional assumption and aiming at a complete study. Notice that it is a particular case of a pp-wave, thus the results of \cite{Leistner-Schliebner13} apply. The local model is well-known, but a lot of global invariants, some of them of dynamical nature, appear.

On the other hand, as in the conformal or projective situation,  if $\afsuris$ is big enough then it should be possible to give a classification. Indeed, in \cite{gromov} Gromov conjectures that rigid geometric structures with large automorphism group are classifiable. Of course, each situation needs its own definition of the word "large". Results of Ferrand in conformal geometry \cite{Ferrand} or Matveev \cite{Matveev}, recently improved in \cite{Zeghib3}, in projective geometry are typical illustrations of this general conjecture.

In the Riemannian case, a classical theorem by Yano \cite{Yano} asserts that $\afsuris$ is discrete. More recently, Zeghib improved this theorem (see \cite{Zeghib2}), so that he could describe the Riemannian affine actions of lattices of semi-simple groups of rank $\geq 2$. It appears  that when  $\afsuris$ is not trivial, the universal cover splits as a product having a non trivial flat factor. In particular, the finite quotients of flat tori are the only compact \emph{Riemannian} $3$-dimensional manifolds such that $\afsuris$ is non trivial (see also Proposition \ref{decomposable}).\medskip

{\bf Our results.} As often, the situation is richer in the Lorentzian case. Indeed, Yano's Theorem does not extend even in dimension $3$. For example, on a compact flat Lorentzian $3$-manifold with nilpotent holonomy (that is not a torus) there exist $1$-dimensional groups of affine transformations that are not isometric. It can be seen from the description in \cite{Auslander} (Theorem 4 and 9) by simply computing the normalizer of the fundamental group in the (classical) affine group (see also the beginning of the proof of Proposition \ref{transfo irra}).

The local expression of metrics admitting a parallel lightlike vector field $X$ is quite simple (see Lemma \ref{local}) but the global problem is a lot richer and it is not clear that a classification is possible (by classification we mean something close to a list of pairwise non isometric metrics such that any metric admitting a parallel lightlike vector field is isometric to one of them). 

The global behaviour of $X$ is also well understood. Zeghib proved in \cite{Zeghib} that it has to preserve a Riemannian metric. It follows from the work of Carri\`ere \cite{Carriere} on Riemannian flows, that, up to finite cover, $X$ is either tangent to the fibre of a circle bundle over a $2$-torus or is a linear vector field on a $3$-torus.

 It allows us to give an expression of any $3$-dimensional compact Lorentzian manifold admitting a parallel lightlike vector field. It turns out that the most interesting case is when $X$ is tangent to the fibre of a circle bundle over a $2$-torus. This torus is naturally endowed with a linear foliation $\overline\FF$ induced by the distribution $X^\perp$ on $M$. When  $\overline\FF$ has a Diophantine or  rational slope then we are able to improve the former description into a complete classification. In fact, we first find, in  Propositions \ref{Diophantine} and \ref{releve},  expressions that are unique up to some finite dimensional group, then we are able to tell when two expressions give isometric metrics. Unfortunately, these criteria do not read nicely; see the discussions distinguishing the isometry classes in Propositions \ref{Diophantine} and \ref{releve}.

Once this description is obtained, we are able to classify the manifolds $(M,g)$ such that the affine group of $g$ is not elementary \ie such that $[\opna{Aff}(g):\opna{Isom}(g)] >2$ and in particular  such that $\afsuris$ is $1$-dimensional  (see Proposition \ref{champs}). All these manifolds admit a periodic parallel lightlike vector field and therefore fibre over the torus.  The situation splits then  in two: either $\overline\FF$ has closed leaves, and the affine group is always non elementary, or  $\overline\FF$ has dense leaves and the metrics with non elementary affine group are given by Proposition \ref{transfo irra}. In both cases, we have an explicit classification of the metrics. 

{\em All these results are gathered in Section \ref{tables} p.\@ \pageref{tables}, notably in two tables constituting two theorems. Table \ref{table_notations} p.\@ \pageref{table_notations} (re)introduces all the notation used in this section, making it autonomously readable.} Theorem \ref{th_table1}, corresponding to Table \ref{table1} p.\@ \pageref{table1}, gives the list of compact indecomposable Lorentzian 3-manifolds having a parallel (lightlike) vector field and points out those with non compact isometry group. Theorem \ref{th_table2}, corresponding to Table \ref{table2} p.\@ \pageref{table2}, lists those manifolds after the form of the quotient $\opna{Aff}(g)/\opna{Isom}(g)$ ---~which may be trivial or isomorphic to $\Z$, $\Z^2$ or $\R$, by different ways each time. As Lorentzian 3-manifolds with $\opna{Aff}(g)/\opna{Isom}(g)$ non elementary \ie $[\opna{Aff}(g):\opna{Isom}(g)] >2$ are, either flat 3-tori, or such indecomposable manifolds, this gives their list. Gathering the results gives also rise 
to some comments, that we give in a series of remarks.

We have moreover obtained a kind of qualitative description of these metrics. Indeed, we proved in Theorem \ref{theo_equivariant}, that the metrics such that $[\opna{Aff}(g):\opna{Isom}(g)] >2$ are equivariant deformations of flat metrics, in the sense that for any $\varphi\in \opna{Aff}(g)$ there exists a path of metrics $g_t$ between $g$ and a flat metric (in fact the quotient by a lattice of a left invariant metric on $\R^3$ or the Heisenberg group) such that  $\varphi\in \opna{Aff}(g_t)$ for all $t$.   Nevertheless, there exist metrics $g$ on parabolic tori such that $[\opna{Aff}(g):\opna{Isom}(g)] >2$ and  $\opna{Aff}(g)$ is not contained in the affine group of any flat metric.

Finally, we generalize the construction of the metrics of this article to show that compact Lorentzian manifolds with large affine group may be easily built in greater dimension, some of them carrying no flat Lorentzian metric.\medskip

{\bf The article is organized as follows:} Section \ref{reducible} shows that when $\afsuris$ is not elementary and $(M,g)$ is not a flat torus then there exists a parallel lightlike vector field $X$ on $(M,g)$; Section \ref{notcompact} deals with the case  where the orbits of $X$ are non periodic and Sections \ref{closed} and \ref{special} with that where all its orbits of are periodic. Section \ref{transfo} describes the manifolds $(M,g)$ such that $[\opna{Aff}(g):\opna{Isom}(g)]>2$ and states and proves Theorem \ref{theo_equivariant}. Section \ref{tables} gathers the results except Theorem \ref{theo_equivariant} and is readable by itself. Section \ref{section_grande_dim} gives examples in dimension greater than three.\medskip

\noindent{\bf General notation.} If $X$ is a lightlike vector field on $M$, $X^\perp$ its orthogonal distribution, and $\pi:M\rightarrow \overline M$ the quotient of $M$ by the flow of $X$ when it exists, $\mathcal X$ denotes the integral foliation of $X$, $\FF$ that of $X^\perp$ and $\overline \FF$ the image of the latter by $\pi$ on $\overline M$.

\section{Reduction to the lightlike case}\label{reducible}

In this section, $(M,g)$ is a compact Lorentzian $3$-manifold and $D$ is its Levi-Civita connection.
\begin{de}
A Lorentzian manifold is called \emph{decomposable} when the representation of its holonomy group preserves an orthogonal decomposition of the tangent bundle, otherwise said the tangent bundle is an orthogonal direct sum of parallel subbundles.  
\end{de}

\begin{rem}\label{rem_un_seul_champ_parallele} If $(M,g)$ admits two linearly independent parallel vector fields $U$ and $V$, it is a flat torus, up to a possible 2-cover. Indeed, if $\opna{Span}(U,V)$ is $g$-non-degenerate, its orthogonal is spanned by a field $W$ with $|g(W,W)|=1$ (here use a possible 2-cover), which is parallel, so $(U,V,W)$ is a parallel frame field. This defines a locally free action of $\R^3$ on $M$, sending the canonical coordinate fields $(\partial_x, \partial_y, \partial_z)$ of $\R^3$ on $(U,V,W)$. Therefore, classically, $M$ is a flat torus, indeed:\medskip

-- this action sends a parallel frame field on parallel frame field, hence sends the (flat) connection on $\R^3$ on the connection $D$ of $M$, hence the latter is also flat,\medskip

-- as $M$ is compact, there is some $\varepsilon>0$ such that for every $m\in \R^3$, this locally free action is free on $\left]m-\varepsilon,m+\varepsilon\right[^3$. Thus this action defines a covering $\R^3\rightarrow M$, the associated action of $\pi_1(M)$ with which necessarily preserves $(\partial_x, \partial_y, \partial_z)$ hence acts by translations. So $M\simeq\R^3/\Z^3$.\medskip

\noindent If $\opna{Span}(U,V)$ is lightlike, one may suppose that $U\in\ker(g_{|\opna{Span}(U,V)})$ and then the field $W$ defined by $g(U,W)=1$, $g(V,W)=g(W,W)=0$ is also parallel; the same conclusion follows.
\end{rem}

\begin{rem}\label{rem_volume_preserve} For any $\varphi\in \opna{Aff}(g)$, $\varphi^*\opna{vol}=\pm\opna{vol}g$. Indeed, as $\varphi$ preserves $D$ and as $\opna{vol}(g)$ is $D$-parallel, there exists a number $l$ such that $\varphi^*\opna{vol}(g)=l\opna{vol}(g)$. But as $M$ is compact, it has a finite volume preserved by $\varphi$ and therefore $|l|=1$. 
\end{rem}

\begin{prop}\label{decomposable} Suppose that $(M,g)$ is decomposable. Then if $[\opna{Aff}(g):\opna{Isom}(g)] >2$, $(M,g)$ is finitely covered by a flat torus.
\end{prop}
{\bf Proof.} Let $\varphi$ be an element of $\opna{Aff}(g)\smallsetminus \opna{Isom}(g)$. As $(M,g)$ is $3$-dimensional and decomposable, there exists a parallel non-lightlike line field; replacing perhaps $M$ by a 2-cover of it, we take on this line a vector field $X$ such that $|g(X,X)|=1$. We denote its flow by $\Phi^t_X$. According to Remark \ref{rem_un_seul_champ_parallele}, we can suppose that $X$ is the unique parallel vector field, up to proportionality, and that therefore $\opna{Span}(X)$ is preserved by $\varphi$. It implies also that the distribution $X^\perp$ is also preserved (it is the only codimension $1$ parallel distribution transverse to the parallel vector field; if there is another one $P$, the orthogonal projection of $X$ onto $P$ would be another parallel vector field). We recall that this distribution is integrable and defines a codimension 
$1$ foliation that we will denote by $\FF$.

The vector field $X$ being parallel it is also Killing. It means that  we can replace $\varphi$ by any $\Phi_X^t\circ\varphi$, therefore we can assume that $\varphi$ preserves any given leaf $\FF_0$ of $\FF$.
 Clearly, $\varphi|_{\FF_0}$ is a diffeomorphism of the Riemannian or Lorentzian surface $(\FF_0, g|_{\FF_0})$ that preserves the induced connection. If $\varphi|_{\FF_0}$ preserves  the metric then, as $\varphi$ preserves the volume form of $g$, it preserves $X$ and then $\varphi$ is isometric contrarily to our assumption. 
Hence $(\FF_0, g|_{\FF_0})$ is a surface that admits a non isometric affine transformation $\varphi$. Let us denote by $S$ the self adjoint automorphism such that $\varphi^\ast g=g(\,\cdot\,,S\,\cdot\,)$. On the one hand, if $S=\lambda\Id$ then $|\lambda|=1$ ---hence $\lambda=-1$ as $\varphi$ is not an isometry. Indeed, the scalar curvature $r$ is the trace of $\Ric$, the endomorphism such that $\ric=g(\,\cdot\,,\Ric\,\cdot\,)$. Now $\varphi^\ast g=\lambda g$ and $\varphi$ preserves $D$ thus$\varphi^\ast\ric=\ric$, consequently $r\circ\varphi=\frac1\lambda r$. But $r$ is bounded as $M$ is compact so $|\lambda|=1$. On the other hand, if $S\not\in\R.{\rm Id}$ then $g|_{\FF_0}$ is flat, by the following fact.\medskip

\noindent{\em Fact.} Let $(N,g)$ be some Lorentzian surface. If it admits a parallel field $S\not\in\R.{\rm Id}$ of self adjoint endomorphisms, then it is flat. Indeed, by standard linear algebra, the couple of matrices of $(g,S)$ is, up to conjugation, of one of the three following types:
\begin{center}
{\small$\displaystyle \left\{\left(\begin{array}{cc}1&0\\0&\pm1\end{array}\right),\left(\begin{array}{cc}\lambda&0\\0&\lambda'\end{array}\right)\right\},\ 
\left\{\left(\begin{array}{cc}0&1\\1&0\end{array}\right),\left(\begin{array}{cc}0&\lambda\\-\lambda&0\end{array}\right)\right\},\ 
\left\{\left(\begin{array}{cc}0&1\\1&0\end{array}\right),\left(\begin{array}{cc}\lambda'&1\\0&\lambda'\end{array}\right)\right\}\ $}
\end{center}
with $0\neq\lambda\neq\lambda'$. In each case, the holonomy algebra is included in $\{A;g(A\,\cdot,\,\cdot\,)=-g(A\,\cdot,\,\cdot)\ \text{and }AS=SA\}=\{0\}$. The fact follows.\medskip

So finally, either $\varphi|_{\FF_0}$ is an anti-isometry and $[\opna{Aff}(g):\opna{Isom}(g)]=2$ ---such a case is realized for example by the product of the Clifton-Pohl torus by $S^1$--- or $[\opna{Aff}(g):\opna{Isom}(g)]>2$ and then $g_{\FF_0}$ is flat, hence also $g$, that reads locally as a product $g|_{\FF_0}\times \dd t^2$.

In the second case, if the vector field $X$ is timelike then the metric obtained simply by changing the sign of $g(X,X)$ is flat and Riemannian therefore $M$ is finitely covered by a torus.
If $X$ is spacelike it follows from the classification of $3$-dimensional Lorentzian flat manifolds (see \cite{Goldman}) that 
either $M$ is finitely covered by a torus or $(M,g)$ has hyperbolic holonomy and is obtained by suspending a hyperbolic automorphism of $\T^2$. In the last case, the flow of $X$ is an Anosov flow 
whose strong stable and unstable direction  are  the lightlike direction of $g|_{\mathcal{F}}$. It follows that $\varphi$ preserves these line fields. 
In that case the leaves of $\FF$ are compact, therefore, by Remark \ref{rem_volume_preserve} the restriction of $\varphi$ to any leaf $\FF_0$  preserves the volume.       The map $\varphi|_{\FF_0}$ preserving the volume and the lightlike direction has to be an isometry. We have a contradiction. \boxfill$\Box$\medskip

Let us remark that if $(M,g)$ is a compact Riemannian $3$-manifold such that $\opna{Aff}(g)\supsetneq \opna{Isom}(g)$ then it is decomposable and thus the proof of Proposition \ref{decomposable} says (because Riemannian anti-isometries do not exist) that $(M,g)$ is finitely covered by a flat Riemannian torus.

\begin{rem}\label{pasdindicedeux} The proof of Proposition \ref{decomposable} shows that the case where $(M,g)$ is decomposable and $\indaf=2$ is very specific; the metric on the universal cover $\wdt M$ is a product $g_1\times g_2$ with $g_1$ Riemannian, 1-dimensional and $g_2$ Lorentzian; $\opna{Aff}(g)/\opna{Isom}(g)$ is spanned by a diffeomorphism preserving $g_1$ and changing the sign of $g_2$. This would demand a specific study; besides, the dynamic phenomena induced by an {\em infinite} group $\opna{Aff}(g)/\opna{Isom}(g)$, in which we are interested here, do not appear in this case. So we do not study it here and the rest of the article deals with the case where $M$ is indecomposable.
\end{rem}

\begin{lem}\label{basic_pres} If $M$ is indecomposable and $\varphi\in\opna{Aff}(g)\smallsetminus\opna{Isom}(g)$, then $\varphi^*g(.,.)=g(E.,.)$ with $E$ parallel, $E=\Id+N$ and $N^2=0$; more precisely, up to a possible 2-cover making $M$ time-orientable, there is a parallel lightlike vector field $X$, canonically associated with $\varphi$, up to sign, and a constant $C\in\R^\ast$ such that $\varphi^\ast g=g+CX^\flat\otimes X^\flat$. Saying it with matrices, if $M$ is moreover orientable, in any frame field $(X,Y,Z)$ in which: 
$$\Mat(g)=\text{\small$\left(\begin{array}{ccc}0&0&1\\0&1&0\\1&0&0\end{array}\right)$}\ \text{(such frames exist), then: }\ \Mat(N)=\text{\small$\left(\begin{array}{ccc}0&0&C\\0&0&0\\0&0&0\end{array}\right)$}.$$
Moreover if $M$, up to a 2-cover, is not a flat torus, there exists $\lambda \in \R^*$ such that $\varphi_\ast X =\lambda X$.
\end{lem}
{\bf Proof.} Some facts proven in \cite{Martinealaplage} may shorten a bit the proof; a nevertheless brief and self-contained is possible, so we give one. Let $P\in\R[X]$ be the minimal polynomial of $E$. If $P=P_1P_2$ with $P_1$ and $P_2$ prime with each orther then $TM=\ker P_1(E){\oplus}\ker P_2(E)$, the sum being orthogonal as $E$ is self adjoint. But this is ruled out as $M$ is indecomposable. So, as $\dim M$ is odd, $P$ is of the form $(X-\alpha)^k$, with $\alpha\in\{\pm1\}$ as $\varphi^*\opna{vol}(g)=\pm\opna{vol}(g)$; actually $\alpha=1$ as $\varphi^\ast g$ must have the same signature as $g$. So $E=\Id+N$ with $N$ nilpotent and non null as $\varphi\not\in\opna{Isom}(g)$. The image of $N$ is clearly parallel and lightlike. Let $h$ be the symmetric $2$-tensor on $\opna{Im} N$ defined by $h(Nu,Nv)=g(u,Nv)$. It is well defined as $u$ is defined up to an element of $\ker N$ and $N$ is symmetric. This tensor is actually non-degenerate \ie a parallel Lorentzian or Riemannian metric on $\opna{Im} N$. Hence, if the distribution $\opna{Im}N$ is $2$-dimensional then it contains two parallel line fields (both non-degenerate or both lightlike for $h$, this does not 
matter), one of them being $\opna{Im}N^2$. But then the other one has to be $g$-non-degenerate, else $\opna{Im}N$ would be $g$-totally isotropic. This contradicts the fact that $(M,g)$ is indecomposable. Thus, $\opna{Im}N$ is an oriented (because $(M,g)$ is time-orientable) parallel lightlike line field. We define $X$ as the vector field, unique up to sign, tangent to $\opna{Im}(N)$ such that $h(X,X)=1$. It is parallel and lightlike and there exists $C>0$ such that $N=[u\mapsto C\,g(u,X)X]$ (this was first proved in \cite{Martinealaplage}) \ie $\varphi^\ast g=g+CX^\flat\otimes X^\flat$. As $M$ is supposed to be orientable, we may extend $X$ as a continuous frame field $(X,Y,Z)$; we may moreover take $Z$ lightlike and $Y\perp\opna{Span}(X,Z)$. Rescaling possibly $Y$ and $Z$, we obtain by construction the announced form for $\Mat(g)$ and $\Mat(N)$.

Finally, as $X$ is parallel, so is $\varphi_\ast X$ and the last assertion follows from Remark \ref{rem_un_seul_champ_parallele}.\boxfill$\Box$\medskip

It follows from Lemma \ref{basic_pres} that investigating the Lorentzian $3$-manifold with $\opna{Aff}(g)\supsetneq\opna{Isom}(g)$ requires to study the indecomposable Lorentzian $3$-manifolds with a parallel, lightlike vector field. It is now our work, shared out among sections \ref{notcompact},  \ref{closed} and \ref{special} after the cases given by Proposition \ref{gz+yc}.\medskip

{\bf From now on, \mathversion{bold}$(M,g)$ is an orientable, indecomposable compact Lorentzian 3-manifold admitting a parallel lightlike vector field $X$ (which is then unique, by Remark \ref{rem_un_seul_champ_parallele}) and that is not, up to a possible 2-cover, a flat torus.\mathversion{normal}}\medskip

The $1$-form $X^\flat$, defined by $X^\flat u=g(X,u)$, is closed and therefore its kernel, which is $X^\perp$, generates a codimension $1$ foliation that we will denote by $\FF$. We will denote by $\mathcal X$ the foliation generated by $X$. The following proposition is the key of the study that follows. We do not know whether its analogue in higher dimension is true.

\begin{rem}\label{rem_phi_preserve_Xbemol} $X^\flat$ is stable by $\varphi$ up to proportionality ---~so in particular, $\varphi^\ast\FF=\FF$. Indeed, $X^\flat$ is the only parallel 1-form with $X$ in its kernel. If $\alpha$ is another one, then $\ker\alpha$ is $g$-non-degenerate, but this is impossible as $M$ is indecomposable.
\end{rem}

\begin{prop}\label{gz+yc}  The  vector field $X$ preserves a frame field and there exists $i\in\{1,2\}$ such that the closure of the orbits of $X$ are $i$-dimensional tori. Moreover, if $M$ is orientable then these tori are the fibres of a fibration over a $(3-i)$-dimensional torus.
\end{prop}
\noindent{\bf Proof.} The vector field $X$ being parallel is also Killing. Non equicontinuous Killing vector fields of compact Lorentzian $3$-manifolds, {\em i.e.\@} fields the flow of which is not relatively compact in $\opna{Diff}(M)$, have been classified by Zeghib in \cite{Zeghib}. It turns out that none of them is parallel and lightlike, and therefore that $X$ preserves a Riemannian  metric $\alpha$ ---built by an averaging process over the closure of $(\Phi^t_X)_t$, which is a compact Lie group in $\opna{Diff}(M)$. We get then more: the vector field $W$ defined by $g(X,W)=1$ and $W\perp_\alpha\FF$ is stable by the flow of $X$, and so is the (unique) field $V$ tangent to $\FF$ and such that $\alpha(V,V)=1$, $\alpha(X,V)=0$ and $(X,V,W)$ is a direct frame field (this last condition is possible as $M$ is orientable). So we get a frame field $(X,V,W)$, preserved by $X$.

It means in particular that the flow of $X$ preserves a frame field of the normal bundle $TM/\XX$.  It follows then from Carri\`ere's article \cite{Carriere} that the closure of the leaves of $\mathcal X$ are tori of constant dimension and that these tori are the fibre of a fibration. The basis of the this fibration is compact and orientable and its Euler number is zero, therefore is is a torus.

Furthermore, if the closure of the leaves  of $\mathcal X$ are $3$-dimensional then $(M,g)$ is a flat torus and therefore it is decomposable.\boxfill$\Box$\medskip

The study naturally splits in two cases according to the dimension of the closure of the orbits of $X$. Before studying them, we state a lemma showing that indecomposable Lorentzian 3-manifolds with a parallel vector field have a unique local invariant, their curvature, which is a scalar. Therefore, once this is set, all the work that follow to get classification results deals only with global questions about the metric $g$, seeking global invariants.

\begin{lem}\label{local} As $X$ is parallel, the curvature tensor of $(M,g)$ is entirely given by the data of $g(R(Z,Y)Z,Y)$ with $(Y,Z)$ such that $Y\perp X$,  $g(Y,Y)=1$ and vol$(X,Y,Z)=1$. The obtained function does not depend on the choice of such a couple; we denote it by $r$. Notice that $R(Y,Z)Y=r X$.
 
There are local coordinates $(x,y,z)$ around each $m\in M$ such that:\medskip

$\bullet$ $\frac\partial{\partial x}=X$ and $\Span\left(\frac\partial{\partial x},\frac\partial{\partial y}\right)=X^\perp$,\medskip

$\bullet$  $\Mat_{(\frac\partial{\partial x},\frac\partial{\partial y}\frac\partial{\partial z})}(g)=\text{\small$\left(\begin{array}{ccc}0&0&1\\0&1&0\\1&0&\mu(y,z)\end{array}\right)$}$ with {\small$\left\{\begin{array}{cr}\text{\bf (i)}&\mu(0,\,\cdot\,)\equiv0\\\text{\bf (ii)}&\frac{\partial\mu}{\partial y}(0,\,\cdot\,)\equiv0\\\end{array}\right.$} .\medskip

\noindent Then $r=\frac{\partial^2\mu}{\partial y^2}$. Therefore, in the neighborhood of $m$, $g$ is locally characterized, up to isometry fixing $m$, by the function $r$. Such coordinates, centered at $m$, are unique once the lightlike vector $\frac{\partial}{\partial z}$ is fixed at $m$. 
\end{lem}
{\bf Proof.} The first assertion is nearly immediate and left to the reader. Take any neighborhood of ${\cal O}$ of $m$ such that $({\cal O},X)\simeq(\R^3, e_1)$. On the quotient $\overline{{\cal O}}$ of $\cal O$ by the orbits of $X$, take a field $\oY$ tangent to $\overline{\FF}$ such that $\gbar(\oY,\oY)=1$. Actually this conditions determines $\oY$ up to its sign. Then take any transversal $\TT$ to $\overline{\FF}$, and $\oZ$ the unique field along $\TT$ such that, for any pull back $Z$ of it on $\cal O$, $g(X,Z)\equiv1$. Then extend $\oZ$ to the whole $\overline{{\cal O}}$ by pushing it by the flow of $\oY$. Then there is a unique pull-back $Y$ of $\oY$ on ${\cal O}$ such that $g(Y,Z)\equiv0$, with $Z$ any pull-back of $\oZ$ ---~this condition being independent of the choice of $\oZ$. So in turn, once chosen any pull-back $Z$ of $\oZ$ along $\TT$, commuting with $X$, there is a unique way to extend it in a pull-back $Z$ commuting with $Y$. Finally, as $X$ is a Killing field for $g$, by construction, $Y$ and $Z$ commute with $X$. Taking $(x,y,z)$ the coordinates dual to $(X,Y,Z)$ gives the wanted 
form for $g$, possibly without Conditions {\bf (i)} and {\bf (ii)} for the moment. Now just above, the choice of the transversal $\TT$, and that of the pull back $Z$ of $\oZ$ along it, were free. So we may choose them such that the orbit of $Z$ through $m$ is a lightlike geodesic.
 This is exactly equivalent to ensuring {\bf (i)} and {\bf (ii)}. The only arbitrary choice is then that of $Z$ at $m$. \boxfill$\Box$\medskip

We will also need the following technical lemma.

\begin{lem}\label{lem_phi_preserve} If $\varphi\in\opna{Aff}(g)\smallsetminus\opna{Isom}(g)$, there is $\lambda\in\R^\ast$ such that $\varphi_\ast X=\lambda X$ and $\varphi^\ast X^\flat=\frac1\lambda X^\flat$. Moreover, if $r$ is the function introduced by Lemma \ref{local}, $r\circ\varphi=\lambda^2 r$.
\end{lem}

{\bf Proof.} Let us denote by $\lambda\in\R^\ast$ the scalar given by Lemma \ref{basic_pres}, such that $\varphi_\ast X=\lambda X$. By Remark \ref{rem_phi_preserve_Xbemol}, $\varphi^*X^\flat=\lambda' X^\flat$. As the only eigenvalue of the endomorphism $E$ appearing in Lemma \ref{basic_pres} is 1, $\lambda'=\frac1{\lambda}$. Now take $(Y,Z)$ a couple of vector fields as in the statement of Lemma \ref{local}, such that $R(Y,Z)Y=rX$. As $\varphi^*X^\flat=\frac1\lambda X^\flat$, $\varphi_\ast Z\equiv\frac1\lambda Z\,[X^\perp]$ and, as $\varphi$ preserves the volume and as $\varphi_\ast X=\lambda X$, necessarily $\varphi_\ast Y\equiv Y\,[X]$. Now, as $\varphi\in\opna{Aff}(g)$, $R$ is $\varphi-equivariant$: $R(\varphi_\ast Y,\varphi_\ast Z)\varphi_\ast Y=(r\circ\varphi^{-1})\varphi_\ast X$, \ie here $R(Y,\frac1\lambda Z)Y=(r\circ\varphi^{-1})\lambda X$, hence $r\circ\varphi=\lambda^2 r$. \boxfill$\Box$

\section{Indecomposable manifolds with a parallel lightlike vector field with non closed orbits}\label{notcompact}
Here $(M,g)$ is a compact indecomposable Lorentzian $3$-manifold which is not, up to a possible 2-cover, a flat torus, endowed with a lightlike parallel vector field $X$ whose 
closure of the orbits defines a codimension $1$ foliation with toral leaves, that we will denote by $\mathcal G$. We do not suppose that $M$ is orientable in this section; it turns out that it will be necessarily a torus, hence orientable. It means that the map $\pi : M\rightarrow  M/\mathcal G\simeq S^1$ is a fibration (we will see that 
$M$ is actually diffeomorphic to a torus).
Again the situation splits  in two cases. Two leaves $\FF_0$ of $\FF$ and $\GG_0$ of $\GG$ that meet are either equal or transverse to each other, as they are both stable by the flow of $X$, whose orbits are dense in $\GG_0$. In the first case, $\FF_0$ is compact, in the other case, $\overline{\FF_0}$ contains $\GG_0$ so $\FF_0$ is not closed. Besides, as $\FF$ is the integral foliation of a closed 1-form, all its leaves are diffeomorphic, see {\em e.g.\@} \cite{candel-conlon}: on the universal cover $\widetilde M$, the leaves are the levels of a submersion $f$ on $\R$. Endow $M$ with any auxiliary Riemannian metric, then the flow of $\frac1{\|\nabla f\|^2}\nabla f$, well-defined on $\widetilde M$ and $\pi_1(M)$-invariant, hence well-defined also on $M$, sends leaf on leaf and hereby provides a diffeomorphism between any two leaves of $\FF$. Hence all the leaves of $\FF$ or none of them is compact. In the first case, $\FF=\GG$; in the second one, $\FF$  and $\mathcal G$ are transverse.  We are going to prove:
\begin{prop}\label{prop_nonclosed}
 If $(M,g)$ is an indecomposable compact Lorentzian $3$-manifold admitting a parallel lightlike vector field with non closed leaves, then:
\begin{enumerate}
\item If the leaves of $\FF$ are compact then $(M,g)$ is isometric to the quotient of $\R^3$ endowed with the metric $\wdt g=2 \Lambda \dd x\dd z+ L^2(z)\dd y^2$, where $\Lambda\in \R$  and $L$ is a $1$-periodic function,  by the group generated by $(x,y,z)\mapsto (x+1,y+\tau,z)$, $ (x,y,z)\mapsto (x,y+1,z)$ and $(x,y,z)\mapsto (x+r_1,y+r_2,z+1)$ where $(\tau,r_1,r_2)\in \R^3$,  $\tau\not\in\Q$. 
\item If the leaves of $\FF$ are non compact then $(M,g)$ is isometric to the quotient of $\R^3$ endowed with the metric $\wdt g=2 \Lambda \dd x\dd z+ L^2 \dd y^2+ \mu(y)\dd z^2$, where $(\Lambda,L)\in \R^2$ and $\mu$ is a $1$-periodic function, by the group generated by $(x,y,z)\mapsto (x+1,y,z+\tau)$, $ (x,y,z)\mapsto (x+r_1,y+1,z+r_2)$ and $(x,y,z)\mapsto (x,y,z+1)$ where $(\tau,r_1,r_2)\in \R^3$,  $\tau\not\in\Q$.
\end{enumerate}

\noindent Moreover, in both cases, $\opna{Aff}(g)=\opna{Isom}(g)=\{\text{translations preserving $L(z)$, respectively $\mu(y)$}\}$.

\end{prop}
{\bf Proof.} 
First, we suppose that the leaves of $\FF$ are compact. Let $\varphi$ be an element of $\opna{Aff}(g)$. Then the scalar $\lambda$ of Lemma \ref{lem_phi_preserve} is $\pm1$. Indeed, here $\FF=\GG$. Let $\pi : M\rightarrow S^1\simeq M/\mathcal G$. As the form $X^\flat$ is invariant by the flow of $X$ and as the leaves of $\mathcal X$ are dense in those of $\mathcal G=\FF$, we know that there exists a $1$-form $\zeta$ on $S^1$ such that $\pi^*\zeta=X^\flat$. 
Moreover, as $\varphi$ preserves $\FF$, there is a diffeomorphism $\overline \varphi$ of $S^1$ such that $\overline \varphi\circ \pi=\pi\circ\varphi$. This means that $\overline\varphi^\ast\zeta= \frac1\lambda\zeta$, so $|\lambda|=1$ as $S^1$ is compact.\medskip

The closure of the flow of $X$ in the group of diffeomorphisms is clearly a $2$-torus. This torus is contained in $\opna{Isom}(g)$ and its orbits are the leaves of the foliation $\FF$.  We choose a $1$-dimensional closed subgroup in it. It defines a $1$-periodic Killing vector field $Y$ tangent to $\FF$. We then take a vector field $Z$ that is lightlike, perpendicular to $Y$ and such that $g(X,Z)$ is constant. The vector fields $X$ and $Y$ being Killing, we see that $[X,Y]=[X,Z]=[Y,Z]=0$ ---~proving that $M$ is a torus. 
We can therefore choose $Z$ so that it induces a $1$-periodic vector field in $M/\FF$. These fields define coordinates $(x,y,z)$  on $\widetilde M$ the universal cover, \ie a diffeomorphism with $\R^3$. 

Replacing eventually  $X$ by some vector $\Lambda X$ proportional to it, we see that the action of the fundamental group of $M$ on $\widetilde M$ in these coordinates has the desired form. Furthermore, there exists a function $L$ such that the metric $g$ reads $g=2\Lambda \dd x\dd z+ L^2(z)\dd y^2$ in these coordinates.

Let $\varphi$ be an element of $\rm{Aff}(g)$.
As $X^\flat=\dd z$, there exists $C\in \R$ such that:
\begin{equation}\label{forme_phi^*g} 
\varphi^*g=g+C X^\flat\otimes X^\flat = 2\Lambda\dd x\dd z+L^2(z)\dd y^2+C\dd z^2.
\end{equation}

\noindent Replacing possibly $\varphi$ by $\varphi^2$ we will assume that $\varphi$ preserves $X$, $X^\flat$ and the volume form (by Lemma \ref{basic_pres} if for some $k\in\Z^\ast$, $\varphi^k$ is an isometry, then $\varphi$ also). Hence, we have $\varphi^*\dd z= \dd z$, 
$\varphi^*\dd y= \dd y +\alpha \dd z$ and $\varphi^*\dd x= \dd x +\beta \dd y+ \gamma \dd z$ for some functions $\alpha$, $\beta$ and $\gamma$. 
We see that $\varphi$ acts by translation on the coordinate $z$. As $\varphi^*g$ must read as in \eqref{forme_phi^*g}, the function $L$ is invariant by this translation.

By the symmetry of second order derivatives, we know that $Y.\alpha=0$ and $X.\beta=X.\alpha=X.\gamma=0$. But any function constant along the orbits of $X$ is also constant along the leaves of $\GG=\FF$, therefore we also have $Y.\beta=Y.\gamma=0$. Using again the symmetry of second order derivatives we get $Z.\beta=0$ \ie $\beta={\rm cst}$. Computing $\varphi^\ast g$ and using \eqref{forme_phi^*g} we get:
\begin{equation}\label{eq_bidule}
 \begin{aligned}
  \Lambda\beta + L^2(z)\alpha(z)=0\\
  2\gamma(z)+L^2\alpha^2(z)=C
 \end{aligned}
\end{equation}
It means that $\wdt \varphi(x,y,z)=(x+\beta y + J(z), y - K(z),z+t)$, for some functions  $J$ and $K$. But, if $\beta\neq 0$,  this map is not in the normalizer of the fundamental group, as $\wdt{\varphi}(x+1,y+\tau,z+t)=(x+\beta y + J(z), y - K(z),z+t)+(1,\tau,0)+(\tau \beta,0,0)$ and  $[(x,y,z)\mapsto(x+\tau\beta,y,z)]$ is not in the fundamental group unless $\beta=0$.  Hence, $\beta=0$. By \eqref{eq_bidule}, $\alpha=0$, $\gamma$ is constant and $\wdt \varphi$ reads:
$$(x,y,z)\mapsto (x+\gamma z +c_1, y+c_2, z+t)$$
where $c_1$, $c_2$ are constants. As $X$ and $Y$ are Killing we can assume that $c_1=c_2=0$. But if $\gamma\neq 0$, this map is not in the normalizer of the fundamental group because: 
$$(x+r_1,y+r_2,z+ 1)\mapsto (x+\gamma z,y,z)+(r_1,r_2,1)+(\gamma ,0,0)$$
and $[(x,y,z)\mapsto(x+\gamma,y,z)]$ is not in the fundamental group unless $\gamma=0$. Hence $\varphi\in \opna{Isom}(g)$.  

\hop

Now, we suppose that $\FF$ has non compact leaves. The closure of the flow of $X$ in the group of diffeomorphisms is still  a $2$-torus contained in $\opna{Isom}(g)$, but now its orbits are the leaves of the foliation $\mathcal G$ that is transverse to $\FF$.
This time, the choice of a closed subgroup in this torus provides us a $1$-periodic Killing vector field $Z$. The quotient $M/\GG$ is 1-dimensional, hence orientable. So we can choose, transversely to $\GG$  a field $Y$ tangent to $\FF$, perpendicular to $Z$ and such that  $g(Y,Y)$ is constant. 
The vector fields $X$ and $Z$ being Killing, we have $[X,Y]=[X,Z]=[Y,Z]=0$ ($M$ is again a torus). We can rescale $Y$ so that it induces a $1$-periodic vector field on $M/\mathcal G$. Take $A=X$ or $A=Y$. Then: $A.(g(X,Z))=g(X,D_AZ)=g(X,D_ZA)=Z.(g(X,A))=0$, as $X$ is parallel. Besides, $Z.(g(X,Z))=({\cal L}_Zg)(X,Z)+g([Z,X],Z)=0$ as $Z$ is Killing and $[Z,X]=0$. Hence $g(X,Z)$ is constant.  We replace $X$ by a vector field proportional to it in order to have a fundamental group of  the desired form.

These vector fields define coordinates $(x,y,z)$ on $\wdt M$ and there exist $(L,\Lambda)\in \R^2$ and a non constant function $\mu$ such that the metric $g$ then reads $g=2\Lambda \dd x\dd z+ L^2 \dd y^2+ \mu(y)\dd z^2$ (if $\mu$ is constant then $Y$ would be parallel, hence $M$ decomposable, contrarily to our assumption.

Take $\varphi$ in $\opna{Aff}(g)$. By Lemma \ref{basic_pres}, 
 $\varphi^*g=g+C X^\flat\otimes X^\flat$. Repeating the reasoning above, but using  the fact that $\varphi$ preserves $\FF$ and $\mathcal G$,  and with the scalar $\lambda$  of Lemma \ref{lem_phi_preserve}, we prove that $\varphi^* \dd x = \lambda \dd x + \beta\dd y +\gamma \dd z$,  $\varphi^* \dd y=\dd y$ and $\varphi^* \dd z = \lambda^{-1} \dd z$, where $\lambda$ is the number given at the beginning of the proof.
Let us show that $|\lambda|=1$.  By Lemma \ref{lem_phi_preserve}, the curvature function $r$ of Lemma \ref{local} satisfies $r\circ\varphi=\lambda^2 r$. Unless $|\lambda|=1$, this compels $r$ to be null at every recurrent point of $\varphi$. Now $M$ is compact and $\varphi$ is volume preserving so almost every point of $M$ is recurrent. But $M$ is not flat (else $Y$ would be parallel, see above), thus $r\not\equiv0$, and $|\lambda|=1$.

Therefore, \eqref{eq_bidule} still holds, up to possible $\pm$ signs, and with this time $L={\rm cst.}$ and $\alpha=0$. It gives $\beta=\gamma=0$ hence $\varphi\in\opna{Isom}(g)$.\boxfill$\Box$

\section{Indecomposable manifolds with a parallel lightlike field with closed orbits}\label{closed}
Now  $(M,g)$ is an orientable, compact indecomposable Lorentzian $3$-manifold endowed with a lightlike parallel vector field with closed orbits $X$. Our goal is to understand, as precisely as possible, the isometry classes of such manifolds. This is done, after introductory remarks, in particular the important Remark \ref{rem_not1}, by Proposition \ref{enforme} and Corollary \ref{classesdisometrie_Xperiodique}; see also the commentaries before and after this corollary.

It follows from Proposition \ref{gz+yc} that the flow of $X$ is periodic (from now on, rescaling $X$, we will assume that it is $1$-periodic) and that $X$ is tangent to the fibre of a fibration $\pi :\, M\rightarrow  \Mbar$, where $\Mbar$ is  diffeomorphic to $\T^2$.
Orientable circle bundles over a torus are well know,  they are diffeomorphic either to a $3$-torus or to a manifold obtained by suspending a parabolic automorphism of $\T^2$. 
\begin{notation}\label{notation_Gamma}
For any $n\in \N$ and $(c_1,c_2)\in [0,1[^2$, we define $\Gamma_{n,c_1,c_2}$ as the group of diffeomorphisms of $\R^3$ generated by the maps 
\begin{eqnarray*}
(x,y,z)&\mapsto& (x+1,y,z)\\
(x,y,z)&\mapsto& (x+c_1,y+1,z)\\
(x,y,z)&\mapsto& (x+ny+c_2,y,z+1)
\end{eqnarray*}
We will denote simply by $\Gamma_n$ the group $\Gamma_{n,0,0}$.
\end{notation}
When $n\neq 0$, the reader will recognize the action of a lattice of $H$, the $3$-dimensional Heisenberg group, on itself. Indeed, this action preserves the moving frame $(\partial_x,\partial_y+nz\partial_x,\partial_z)$ and the Lie algebra spanned by these three vector fields is the Lie algebra of $H$. It is not a surprise as  orientable, non trivial circle bundles over the torus are quotient of the Heisenberg group (what is called a Nil-manifold). Besides, if $\widetilde g$ is a left invariant metric on the Heisenberg group such that the center of $H$ is lightlike, then it is flat and it admits a parallel lightlike vector field. We can translate this situation in our vocabulary. Let $\widetilde g$ be the metric on $\R^3$ whose matrix in the moving frame $(\partial_x,\partial_y+nz\partial_x+\theta \partial_z,\partial_z)$ is  
\begin{eqnarray}\label{exemples plats}
\begin{pmatrix}
0       & 0        &  \Lambda\\
0       &  L^2   &  N \\
\Lambda & N & M
\end{pmatrix},
\end{eqnarray}
where $\theta$, $\Lambda$, $L$, $M$ and $N$ are constants. The metric $\widetilde g$ induces a flat metric $g$ on $\R^3/\Gamma_{n,c_1,c_2}$ and the image of $\partial_x$ is a well defined lightlike and parallel vector field. In order to obtain non flat examples, we can replace the constants $L$, $M$, and $N$ by functions of $y$ and $z$ invariant by the action of $\Gamma_{n,c_1,c_2}$ \ie (1,1)-biperiodic. After such a change, the field $\partial_x$ remains parallel.

\begin{prop}\label{fibres}
The map $\pi : \R^3/\Gamma_{n,c_1,c_2} \rightarrow \R^2/\Z^2$ with $\pi([x,y,z])=([y,z])$ is a circle bundle whose Euler number is $n$. Hence,  any orientable circle bundle over the torus is isomorphic to a $\pi_{n,c_1,c_2}$ (actually $c_1$ and $c_2$ do not play any role here).

Let $\widetilde g$  be the metric on $\R^3$ whose matrix in the moving frame $(\partial_x,\partial_y+nz\partial_x+\theta \partial_z,\partial_z)$ is  
\begin{eqnarray}\label{description bas}
\begin{pmatrix}
0       & 0        &  \Lambda\\
0       & L^2(y,z)   &  \nu(y,z)\\
\Lambda & \nu(y,z) & \mu(y,z)
 \end{pmatrix}, 
\end{eqnarray}
where $L$, $\mu$, $\nu$ are biperiodic functions. 
The Lorentzian manifold $(\R^3/\Gamma_{n,c_1,c_2}, g)$, where $g$ is the metric induced by $\widetilde g$, admits a parallel lightlike vector field with closed orbits.
\end{prop}
\begin{rem}
 When $n\neq 0$, the metric $\wdt g$ defined in Proposition \ref{fibres} can be seen on as a metric on the Heisenberg group expressed in a left invariant moving frame.
\end{rem}

\begin{rem-notation}\label{rem_not1}
Any object invariant, or invariant $\mod X$, by the flow of $X$ defines an object on $\Mbar$. Here are some of them:
\begin{itemize}
\item the connection $D$ gives a connection $\overline D$ on $\Mbar$;
\item the $D$-parallel $1$-form $X^\flat$ gives a $\oD$-parallel $1$-form $\zeta$;
\item the foliation $\FF$ tangent to $\ker X^\flat$ gives a foliation $\overline{\FF}$ tangent to $\ker \zeta$;
\item the $D$-parallel $2$-form $i(X)\operatorname{vol}(g)$ gives a $\oD$-parallel  volume form $\overline \omega$;
\item the restriction of $g$ to $\FF$ gives a Riemannian metric on $\overline{\FF}$, which we denote by $\overline g$. Of course $\oD \overline g=0$.
\item any vector field $Y$ tangent to  $\ker X^\flat$ and such that $g(Y,Y)$ is constant is projected onto a $\overline D$-parallel vector field on $\Mbar$
\end{itemize}
 \end{rem-notation}

Let us see how the connections $\oD$ on the 2-torus, associated with metrics appearing in Proposition \ref{fibres} look like. We denote by $\oY$ and $\oZ$ the projection on $\R^2/\Z^2$ of $\partial_y+nz\partial_x+\theta \partial_z$ and $\partial_z$.  By a direct computation we see that:
\begin{equation}\label{eq_connection}\oD_{\oY}\oY= \frac{\dd L}{L}\, \oY \qquad \text{and}\qquad \oD_{\oZ}\oZ=\frac{\oZ.\nu-\frac{1}{2} \oY.\mu+ n\Lambda}{L^2}\,\oY.\end{equation}

The foliation $\FF$ being given by a closed $1$-form all its leaves are diffeomorphic. Actually there are two possibilities: either they are all diffeomorphic to a torus or they are all diffeomorphic to a cylinder. Moreover in the latter case they are all dense in $M$.
\begin{de}
A (torsion free) connection $\oD$ on $\T^2$ will be called unipotent if its holonomy group is unipotent i.e.\ if  it admits  a non trivial parallel vector field and a parallel volume form.
\end{de}
Hence, the connection $\oD$ induced by $D$ is unipotent. Actually, all the items of Remark \ref{rem_not1} are the consequence of this fact.
\begin{prop}\label{GaussBonnet}
Let $\oY$ and $\oZ$ be vector fields on $\Mbar$ such that $\oY$ is $\oD$-parallel ---~thus tangent to $\overline{\FF}$~--- and $\zeta(\oZ)=\overline g(\oY,\oY)=1$. The field $\oY$ is unique up to multiplication by a scalar $\lambda$.
Then the curvature of $\oD$ is encoded by a single function $r:\Mbar\rightarrow \R$ defined by
$$\oR(\oZ,\oY)\oZ=r\oY,$$
where $\overline R$ is the curvature of $\oD$.
Replacing $\oY$ by $\lambda\oY$ turns $r$ into $\frac1{\lambda^2}r$ so $r$ is well-defined up to a positive factor. The function  $r$ is the only local invariant of $(M,g)$, in particular $r=0$ if and only if  $g$ is flat. Moreover this function satisfies a Gauss-Bonnet equality \ie
$$\int_{\Mbar} r\ombar=0,$$
for any parallel volume form $\ombar$.
\end{prop}
{\bf Proof.} The assertion about the curvature of $\oD$ is obtained by direct computation and the fact that $r$ is only local invariant is given by Lemma \ref{local}. Let us prove the Gauss-Bonnet relation. Let $\oY$ and $\oZ$ be two vector fields on $\Mbar$ such that $\oY$ is  tangent to $\overline{\FF}$ and $\overline g(\oY,\oY)=1$  and such that  $\zeta(\oZ)=1$. Let $\kappa_1$ and $\kappa_2$ be the functions on $\T^2$ defined by $\oD_{\oY}\oZ=\kappa_1 Y$ and $\oD_{\oZ}\oZ=\kappa_2 Y$. A direct computation gives: 
$$r=\oZ.\kappa_1 -\oY.\kappa_2 -\kappa_1^2.$$
Let $\xi$ be the form such that $\xi(\oZ)=0$ and $\xi(\oY)=1$. The parallel form $\overline \omega$ (see Notation \ref{rem_not1}) reads $\ombar=\xi\wedge \zeta$. Moreover, $\dd (\kappa_1 \xi+\kappa_2\zeta)(\oY,\oZ)=r$, therefore   $\dd (\kappa_1 \xi+\kappa_2\zeta)=r\ombar$. It follows by Stokes theorem that $\int_{\Mbar}r\ombar=0$.
\boxfill$\Box$
\begin{prop}\label{fibre}
The map that associates with a compact orientable $3$-dimensional Lorentzian manifold $M$ endowed with a parallel lightlike periodic vector field $X$, a unipotent connection on $\T^2$, identified with $M/{\cal X}$, is onto. 
\end{prop}
{\bf Proof.}  Let $\oD$ be a unipotent connection on $\T^2$. Let $\oY_0$ be a non trivial $\oD$-parallel vector field tangent to $\overline \FF$ and $\overline \omega$ be a parallel volume form. The form $\zeta':=i_{\oY_0}\overline{ \omega}$ is also parallel and therefore closed. By definition $\zeta'(\oY_0)=0$.

We choose coordinates $(v,w):\T^2\rightarrow \R^2/\Z^2$ such that $\zeta'$ has constant coefficients in the basis $(\dd v, \dd w)$ and such that $\int \dd v\wedge \dd w=1$. 
Let $\theta$ be a real number such that $\zeta'(\partial_v+\theta\partial_w)=0$. We will denote now by $\oY$ the vector field $\partial_v+\theta\partial_w$ and by $\oZ$ the vector field $\partial_w$.  There exists a function $a: \T^2\rightarrow \R$ such that $\oY=e^a \oY_0$ and therefore 
$$\oD_.(\oY)= da(.) \oY.$$
The function $a$ actually depends on the choice of $\oY_0$, but it is canonical up to an additive constant.
Furthermore, $\zeta'(\oD_.\oZ)=\oD_.(\zeta'(\oZ))+(\oD_.\zeta')(\oZ)=\oD_.(1)+0=0$. It means that $\oD_{\oZ}\oZ$ is proportionnal to $\oY$ so there exists a function $b:\T^2\rightarrow \R$ such that:
$$\oD_{\oZ}\oZ=\frac{b}{e^{2a}}(\oY).$$
We set 
$$C:=\int_{\T^2}b\, \dd v\wedge\dd w.$$
The form $\kappa=(b-C)\, \dd v\wedge\dd w$ is therefore exact. Let $\delta$ be a $1$-form such that $\dd\delta=\kappa$. Take  $\nu$ and $\mu$ the functions such that  $\delta=-\nu \dd v - \mu/2 (-\theta \dd v+ \dd w).$  
As $(\dd v,-\theta \dd v+ \dd w)$ is dual to $(\oY,\oZ)$:
$$b=\oZ.\nu-\oY.(\mu/2)+C.$$
Now we choose $n\in \N$ and $\Lambda \in \R^*$ such that $n\Lambda=C$. We define now the metric $\widetilde g$ on $\R^3$ given in the moving frame $(\widetilde X,\widetilde Y,\widetilde Z)=(\partial_x, \partial_v+nw\partial_x+\theta \partial_w,\partial_w)$ by  
\begin{eqnarray}\label{description basique}
\begin{pmatrix}
0       & 0        &  \Lambda\\
0       & e^{2a(v,w)}   &  \nu(v,w)\\
\Lambda & \nu(v,w) & \mu(v,w)
  \end{pmatrix}.
\end{eqnarray}
For any $(c_1,c_2)\in \T^2$, this metric 
 induces a metric on $\R^3/\Gamma_{n,c_1,c_2}$ which has the desired properties by \eqref{eq_connection}.\boxfill$\Box$

\begin{rem}
 The constant $C$ appearing during the proof of Proposition \ref{fibre} is null if and only if $\oD$ is the image of a Lorentzian connection on $\T^3$, else it is the image of a connection on a parabolic torus.
\end{rem}
\begin{prop}\label{enforme}
If $(M,g)$ is an orientable Lorentzian $3$-dimensional manifold endowed with a parallel lightlike vector field with closed orbits $X$, then it is isometric to $(\R^3/\Gamma_n,g)$ where $g$ is the metric induced by a metric $\widetilde g$ on $\R^3$ that reads like \eqref{description bas}. 
\end{prop}
{\bf Proof.} Let $(M,g)$ be a Lorentzian metric admitting a parallel lightlike vector field with closed orbits that induces the connection $\oD$ on $\T^2$. We endow $\T^2$ with the coordinates and the vector fields $\oY$ and $\oZ$ defined  at the beginning of the proof of Proposition \ref{fibre} (in particular it gives a real number $\theta$). As orientable circle bundles over the torus are determined by their Euler number, using the first sentence of Proposition \ref{fibres}, we may assume that $M=\R^3/\Gamma_n$ with $X$ tangent to the fibres of $\pi:\R^3/\Gamma_n\rightarrow\R^2/\Z^2$, \ie $X=\partial_x$. Now we must prove that the identification $M=\R^3/\Gamma_n$ may be chosen such that the metric $g$ reads as announced. To do this, first we build a morphism from $\R^3/\Gamma_{n,c_1,c_2}$ to $M$ (identified with $\R^3/\Gamma_n$) having particular properties with respect to the projection $\R^3/\Gamma_{n,c_1,c_2}\rightarrow{\mathbb T}^2$, then we show that we can manage to choose $c_1=c_2=0$.
This is done through Lemmas \ref{relevons}, \ref{hophophop} and \ref{cohomology}, which we state autonomously as they will also be useful farther.
\begin{lem}\label{relevons}
Let $n\in \N$ and $\pi:\R^3/\Gamma_{n,c_1,c_2} \rightarrow \R^2/\Z^2$ be the fibration defined by $\pi([x,y,z])=([y,z])$. Let $(\overline V,\overline  W)$ be a commuting (\ie $[\overline V,\overline  W]=0$) moving frame on $\R^2/\Z^2$. If the volume form on $\R^2/\Z^2$ induced by the moving frame $(\overline V,\overline  W)$ has total volume $1$ then there exists vector fields $V$ and $W$ on $\R^3/\Gamma_{n,c_1,c_2}$ sent respectively by $\pi$ on $\overline V$ and $\overline  W$ and such that $[V,W]=n\partial_x$, $[V,\partial_x]=[W,\partial_x]=0$.

Other fields $(V',W')$ have the same property if and only if there is a closed $1$-form $\beta$ on $\R^2/\Z^2$ such that $V'=V+\beta(\overline V)\partial_x$ and $W'=W+\beta(\overline W)\partial_x$.
\end{lem}
{\bf Proof.}
Take $(a,b,c,d)$ the functions such that $\Vbar=a\partial_y+b\partial_z$ and $\Wbar=c\partial_y+d\partial_z$. We denote by   $V_0$ the image in $\R^3/\Gamma_{n}$ of  $\partial_y+nz\partial_x$ and $W_0$ that of $\partial_z$. Any fields $(V,W)$ commuting with $\partial_x$ that are  pushed by $\pi$ on $(\Vbar,\Wbar)$ read $(a V_0  +b W_0 +\alpha(\Vbar)\partial_x,c V_0 +d W_0 +\alpha(\Wbar)\partial_x)$ with $\alpha$ some $1$-form on $\R^2/\Z^2$.
Then, as $[\Vbar,\Wbar]=0$, $$[V,W]=(ad-bc)[V_0,W_0]+\dd\alpha(\Vbar,\Wbar)\partial_x=\dd y\wedge\dd z(\Vbar,\Wbar)n\partial_x+\dd\alpha(\Vbar,\Wbar)\partial_x.$$
 Now, we want $[V,W]=n\partial_x=\dd v\wedge\dd w(\Vbar,\Wbar)n\partial_x$, where $\dd v\wedge\dd w$ is the dual volume form of the frame $(V,W)$.
 Hence, $(V,W)$ is as wanted if and only if $\dd\alpha=n(\dd v\wedge\dd w-\dd y\wedge\dd z)$. By assumption, $\int_{\T^2}\dd v\wedge\dd w=1=\int_{\T^2}\dd y\wedge\dd z$, therefore $\dd v\wedge\dd w-\dd y\wedge\dd z$ is exact, thus a form $\alpha$ as wanted exists, and is defined up to the addition of a closed $1$-form $\beta$ as announced in the lemma.
 \boxfill$\Box$\hip

According to Lemma \ref{relevons}, there exists lifts of $Y$ and $Z$ of $\oY$ and $\oZ$ such that $[X,Y]=[X,Z]=0$ and $[Y,Z]=nX$. 
\begin{lem}\label{hophophop}
Let $\pi : M\rightarrow \T^2$ be a  circle bundle, $U$ a $1$-periodic vector field tangent to the fibre of $\pi$ and $(U,V,W)$ a moving frame such that $[U,V]=[U,W]=0$, $[U,W]=nU$ for some $n\in \N$. 
If the projections $\Vbar$ and $\Wbar$ of $V$ and $W$ on $\T^2$ are coordinate vector fields on $\T^2$ then there exists a diffeomorphism between $\widetilde M$, the universal cover of $M$, and $\R^3$ such that $\wdt U=\partial_x$, $\wdt V=\partial_y+nz\partial_x$, $\wdt W=\partial_z$, conjugating the fundamental group of $M$ with $\Gamma_{n,c_1,c_2}$ for some $c_1$, $c_2$ in $\left[0,1\right[$ and sending some preimage of the origin $0\in\T^2$ on $0\in\R^3$.
\end{lem}
{\bf Proof.} The proof is written for the case $n\neq 0$; the case $n=0$ is obtained by replacing the Heisenberg group by $\R^3$.
The vector fields $\wdt U$, $\wdt V$ and $\wdt W$  satisfy $[\wdt U,\wdt V]=[\wdt U,\wdt W]=0$ $[\wdt V,\wdt W]=\wdt U$.
Hence, they generate a finite dimensional Lie algebra isomorphic to the Heisenberg one. As these vector fields are complete, there exists  a diffeomorphism between $\wdt M$ and the Heisenberg group $H$ sending $\wdt U$, $\wdt V$ and $\wdt W$ on left invariant vector fields (see \cite{Palais} by Palais).  It follows that  we can assume that the diffeomorphism is between $\wdt M$ and $\R^3$ and that it sends $\wdt U$ on $\partial_x$, $\wdt V$ on $\partial_y+nz\partial_x$ and $\wdt W$ on $\partial_z$. We consider them now as left invariant vector fields of $H$. 

It follows from our assumption that the cube $[0,1]^3$ is a fundamental domain for the action of the fundamental group of $M$ on $\wdt M$. 
As, moreover, the elements of the fundamental group of $M$ preserve the vector fields $\partial_x$, $\partial_y+nz\partial_x$ and  $\partial_z$, it follows that these elements are translations by elements of $H$ so the  fundamental group is a lattice of $H$ that has to be a $\Gamma_{n,c_1,c_2}$. Be carefull that the choice of the moving frame identifies $\R^3$ with $H$ seen as {\footnotesize $\left\{\left(\begin{array}{ccc}1&z&x\\0&1&ny\\0&0&1\end{array}\right) \right\}$}.
\boxfill$\Box$ \medskip

\begin{rem} The constants $c_1$ and $c_2$ appearing in Lemma \ref{hophophop} are the numbers such that the orbits of $V+c_1U$ and $W+c_2U$ through any point of $\pi^{-1}(0)$ are $1$-periodic.
\end{rem}

We can apply Lemma \ref{hophophop} to the moving frame $(X,Y-\theta Z, Z)$. Thus we have an identification of $\wdt M$ with $\R^3$ such that $\widetilde X=\partial_x$, $\widetilde Y=\partial_y+nz\partial_x+\theta \partial_z$ and $\widetilde Z=\partial_z$. The last thing to prove is the fact that we can choose the lift of $\oY$ and $\oZ$ in such a way that $c_1=c_2=0$. It is proven in the following lemma.

\begin{lem}\label{cohomology}
In Lemma \ref{hophophop}, suppose that $V'$ and $W'$ are other lifts of the fields $\Vbar$ and $\Wbar$. Then there is a closed $1$-form $\beta$ such that $V'=V+\beta(\Vbar)U$ and $W'=W+\beta(\Wbar)U$, see Lemma \ref{relevons}.
The values of $c_1$ and $c_2$ given by the application of Lemma \ref{hophophop} to $V'$ and $W'$ depend only on the cohomology class of $\beta$. Furthermore, it is possible to choose $\beta$ such that $c_1=c_2=0$.
\end{lem}
{\bf Proof.} We lift everything once again to $\widetilde M$. According to Lemma \ref{hophophop} there exists an identification between $\widetilde M$ and $\R^3$ such that: 
$(\widetilde{U},\widetilde{V},\widetilde{W}) =(\partial_x, \partial_y+nz\partial_x+\theta \partial_z, \partial_z)$. 
Therefore, 
$\widetilde V'=\partial_y+nz\partial_x+\theta \partial_z +\widetilde\beta(\partial_y+\theta \partial_z)\partial_x $ and 
$\widetilde W'=\partial_z+ \widetilde\beta(\partial_z)\partial_x$. Let $f$ be a function such that $\dd f=\widetilde \beta$, of course $\partial_x f=0$.

We define now a map $F:\R^3\rightarrow \R^3$  by $F(x,y,z)=(x+f(y,z),y,z)$. 
We see that  $F_*(\widetilde U,\widetilde V,\widetilde W)= (\widetilde U,\widetilde V',\widetilde W')$. It means that $F$ is a diffeomorphism between the identifications of $\widetilde M$ with $\R^3$ associated with the moving frames $(\widetilde U,\widetilde V,\widetilde W)$ and $(\widetilde U,\widetilde V',\widetilde W')$.
Moreover, $F$ normalizes the fundamental group of $M$ (that reads as a $\Gamma_{n,c_1,c_2}$) if and only if for any $(p,q)\in \Z^2,\  f(y+p,z+q)-f(y,z)\in \Z$ \ie if the cohomology class of $\beta$ is integer.

If we take now $\beta=-c_1\dd y - c_2 \dd z$ then  $F(x,y,z)=(x-c_1y -c_2 z,y,z)$. This time $F \Gamma_{n,c_1,c_2} F^{-1}=\Gamma_n$. Otherwise said, lifting  $V'$ and $W'$ instead of $V$, $W$ on $\wdt M$ identifies the fundamental group of $M$ with $\Gamma_n$ and therefore $M$ with $\R^3/\Gamma_n$.
\boxfill$\Box$\medskip

This completes the proof of Proposition \ref{enforme}.\hfill$\Box$\medskip

We know now all the manifolds $(M,g)$ admitting a lightlike parallel vector field $X$ with closed leaves. The question is now to determine how many isometric classes of metrics there are. The quotient $M/{\mathcal X}\simeq{\mathbb T}^2$ is endowed with a quotient connection $\overline D$ with unipotent holonomy, see Remark \ref{rem_not1}. Corollary \ref{classesdisometrie_Xperiodique} shows that this quotient contains almost all the information: such a manifold $(M,g)$ is determined, up to isometry, by the isomorphism class of a torus ${\mathbb T}^2$ with a unipotent connection $\overline D$ (we remind that any such torus is obtained as a quotient $((M,g)/{\mathcal X})$, see Proposition \ref{fibre}) and a quadruple of real numbers. In turn, tori $({\mathbb T}^2,\overline D)$ depend of a reparametrization of a linear flow on ${\mathbb T}^2=\R^2/\Z^2$, which is a complicated object, see below Corollary \ref{classesdisometrie_Xperiodique}.

\begin{cor}\label{classesdisometrie_Xperiodique}
Let $g$ be the metric on $M=\R^3/\Gamma_n$ such that the matrix of $\widetilde g$ in the moving frame $(\partial_x,\partial_y+nz\partial_x+\theta \partial_z,\partial_z)$ is:
\begin{equation}\label{desc_bas}
\begin{pmatrix}
0       & 0        &  \Lambda\\
0       & L^2(y,z)   &  \nu(y,z)\\
\Lambda & \nu(y,z) & \mu(y,z)
 \end{pmatrix}. 
\end{equation}
Let $g'$ be another metric on $M$ that admits a lightlike parallel vector field with closed orbits. 
If the connections $\oD$ and $\oD'$ induced by $g$ and $g'$ on the quotients by these orbits (both diffeomorphic to $\T^2$) are isomorphic then there exists a metric $g''$ isometric to $g'$  and $(k_1,k_2,k_3,k_4)\in \R^*\times \R^3$ such that the matrix of  $\widetilde g''$ in the moving frame $(\partial_x,\partial_y+nz\partial_x+\theta \partial_z,\partial_z)$ is:
\begin{equation*}
\begin{pmatrix}
0       & 0        &  k_1\Lambda\\
0       & e^{k_2}L^2(y,z)   &  \nu(y,z)+k_3\\
k_1\Lambda & \nu(y,z)+k_3 & \mu(y,z)+k_4
 \end{pmatrix}. 
\end{equation*}
Moreover, if $n\neq 0$ then $k_1=1$.
\end{cor}
{\bf Proof.} We repeat the constructions from the proofs of Propositions \ref{fibre} and \ref{enforme} for the metrics $g$ and $g'$ starting with the same connection on $\T^2$. Indeed, we choose coordinates $(v,w)$ on $\T^2$, a parallel vector field $\oY_0$ and a moving frame $(\oY,\oZ)$ as in the proof of Propositon \ref{fibre}. These choices give the functions $a$, $b$, the constant $C$.
 Choosing the lift $(Y,Z)$  of $(\oY,\oZ)$ giving the expression \eqref{desc_bas} of $g$, we see that the expression of $g'$ has to be: 
\begin{equation*}
\begin{pmatrix}
0       & 0        &  k_1\Lambda\\
0       & e^{k_2}L^2(y,z)   &  \nu'(y,z)\\
k_1\Lambda & \nu'(y,z) & \mu'(y,z)
 \end{pmatrix}. 
\end{equation*}
for some constants $k_1$ and $k_2$. Using \eqref{eq_connection} we see that the constant $C$ determines $\Lambda$ when $n\neq 0$ and therefore $k_1=1$ in this case. As we saw the functions $\mu$ and $\nu$ (resp. $\mu'$ and $\nu'$) define a form $\delta$ (resp. $\delta'$) such that $\dd \delta=\kappa$ (resp. $\dd \delta'=\kappa$). It means that $\delta-\delta'$ is closed. Hence, there exists an exact $1$-form $\beta$ such that $\delta-\delta'-\beta$ has constant coefficient. We replace now $(Y,Z)$ by $(Y',Z')=(Y+\beta(\oY)\partial_x,Z+\beta(\oZ)\partial_x)$.  The expression of $g'$ in the moving frame $(\partial_x,Y',Z')$ has now the right form. The corollary then follows from Lemmas \ref{hophophop} and \ref{cohomology}.
\boxfill$\Box$\medskip

Unfortunately, determining the isomorphism classes of unipotent connections is not an easy task. It contains as a subproblem, and  first difficulty,  the study of the reparameterizations of linear flows. 
Indeed, let $\oD_1$ and  $\oD_2$ be two  unipotent connections on $\T^2$. If there exists a diffeomorphism sending $\oD_1$ on $\oD_2$ then it sends the parallel vector fields of $\oD_1$ on the set of parallel vector fields of $\oD_2$. Namely, these fields are equal, up to a constant factor, to some field of constant slope in ${\mathbb T}^2\simeq\R^2/\Z^2$, {\em i.e.\@} to a reparametrization of some linear (= constant) field of $\R^2/\Z^2$. We remind that in our case,  on $\overline M=M/\FF$, a $\overline D$-parallel field if the projection $\overline Y$ of a field $Y$ of constant norm tangent to $\XX^\perp$; its integral foliation, with constant slope, is the foliation $\overline{\FF}$.

This means that being able to tell when two unipotent connections are isomorphic implies being able to tell when two reparametrizations of a linear flow are isomorphic. But it comes out of the  works on reparametrization of linear flows with  Liouvillean slope (see for example \cite{Bassam}) that there is no hope for a classification. However, in the other cases, \ie when the slope is rational or Diophantine this problem is simpler and we are able to improve our description.
\section{The special cases}\label{special}
\subsection{When the slope is Diophantine.}

The foliation $\overline{\FF}$ is given by a closed $1$-form, hence there exists coordinates on $\T^2$ such that it is a linear foliation with a given slope.  Of course this number is only defined modulo the action of $\mathrm{SL}(2,\Z)$ but the property of being Diophantine is invariant by this action. Recall that a real number $\theta$ is said to be Diophantine if it is badly approximated by rationals \ie there exist some positive numbers $k$ and $c$ such that for any $(p,q)\in \Z\times \N^*$
$$|\theta p -q|> \frac{c}{q^k}.$$

When the foliation $\overline {\FF}$ has a Diophantine slope then we are able to improve Proposition \ref{enforme} and describe pairwise non isometric models.
\begin{prop}\label{Diophantine}
Let $(M,g)$ be an orientable Lorentzian $3$-dimensional manifold endowed with a parallel lightlike vector field with closed leaves $X$. If the foliation $\overline \FF$ has a Diophantine slope then it is isometric to $(\R^3/\Gamma_n,g')$ where $g'$ is the metric induced by a metric $\widetilde g'$ on $\R^3$ that reads:
$$
\begin{pmatrix}
  0       & 0         & \Lambda\\
  0       & L^2       & k       \\
  \Lambda & k        & \mu(y,z)
  \end{pmatrix}
$$
in the moving frame $(\partial_x,\partial_y+nz\partial_x+\theta \partial_z,\partial_z)$, where $\Lambda\in\R^\ast$, $L\in\R_+^\ast$, $k\in [0,|\Lambda|[$ and $\mu$ is a $(1,1)$-biperiodic function such that $z\mapsto \int_0^1 \mu(y,z)\dd y$ is a constant belonging to $[0,2|\Lambda|[$. Moreover, if $n\neq 0$ we can assume that $k=0$ and $\forall z\in\R,\,\int_0^1 \mu(y,z)\dd y=0$.

Modulo the action of $\Gamma_n$, such coordinates $(x,y,z)$ on $\wdt M$ are unique up to the action of:\medskip

-- $(\Z/2\Z\times\opna{GL}_2(\Z))\ltimes\R^3/\Z^3$ in case $n=0$,\medskip

-- $\opna{GL}_2(\Z)\ltimes (\R/\Z\times(\Z/n\Z)^2)$ in case $n\neq0$,\medskip

\noindent where the normal factor acts by translation $(x,y,z)\mapsto(x+x_0,y+y_0,z+z_0)$ with $(x_0,y_0,z_0)\in\R^3/\Gamma_0=\R^3/\Z^3$ if $n=0$ and $(x_0,y_0,z_0)\in\R/\Z\times(\frac1n\Z/\Z)^2$ else.  The group   $\opna{GL}_2(\Z)$ acts naturally by $\text{\footnotesize$\left(\begin{array}{cc}a&b\\c&d\end{array}\right)$}.(y,z)=(ay+bz,cy+dz)$ on $\pi(M)\simeq\R^2/\Z^2$, and $\Z/2\Z$ changes the sign of $x$. Thus, $(\theta,\Lambda,L,k,\mu)$ characterizes $g$ up to this action, and the affine structure on $\pi(M)\simeq\T^2$ is canonical.

The translations preserve $(\theta,\Lambda,L,k)$ and act by natural composition on the function $\mu$. $\opna{GL}_2(\Z)$ acts as $\text{\footnotesize$\left(\begin{array}{cc}a&b\\c&d\end{array}\right)$}.\theta=\frac{a+b\theta}{c+d\theta}$ on $\theta$; see its action on $(\Lambda,L,k,\mu)$ at the end of the proof.
\end{prop}

\begin{rem}One could also impose $\Lambda>0$; the group action would become a bit smaller, but a bit less natural.
\end{rem}

{\bf Proof.} Let $\oU$ be a non trivial parallel vector field on $\T^2$ that is tangent to $\overline{\FF}$. We claim that, as  $\overline{\FF}$ has a Diophantine slope,  there exists coordinates $(y,z)$ on $\T^2$ such that $\oU=\partial_y+\theta\partial_z$. It follows from the following Theorem by Kolmogorov (see \cite{Kol}, the following statement is coming from \cite{Katok-Robinson} Theorem 3.5)
\begin{te}[Kolmogorov] If $\theta$ is Diophantine then for every smooth function $h$ on $S^1$ there exists a smooth function $\psi$ such that
$$h(z)-\int_{S^1} h d \mu = \psi (R_\theta z) - \psi (z)$$
where $\mu$ is the Lebesgue measure on $S^1$ and $R_\theta$ the rotation of angle $\theta$.
Furthermore, Diophantine numbers are the only ones for which this property
holds.
\end{te}
It is well known (see \cite{Bassam} for example) that the flow of $\oU$ can be seen as the quotient flow of the flow by horizontal translation on $\R\times S^1$ under the identification $(t,z )\sim (t+h(z),R_\theta z)$ where $h$ is a positive smooth function from $S^1$ to $\R$ called the ceiling function. 
We put the coordinates $(y,z)=(t+\psi(z),z)$, where $\psi$ is the function given by Kolmogorov's Theorem,  on $S^1\times \R$. The function $h$ is then replaced by a constant, therefore $\oU$ is linear in the coordinates induced by $(y,z)$ on the torus. It proves our claim. Moreover, the $1$-form $\zeta$ induced on $\T^2$ by $X^\flat$ has automatically constant coefficients in the coordinates $(y,z)$.

We repeat now the proof of Proposition \ref{enforme}, but we start with the coordinates $(y,z)$ found above. In particular we choose  lifts $\wdt U$ and $\wdt Z$ of $\oU$ and $\partial_z$ to $\wdt M$.
We obtain  coordinates $(x,y,z)$ on $\wdt M$ such that the fundamental group is $\Gamma_n$ and  that the metric reads 
$$
\begin{pmatrix}
  0       & 0         & \Lambda\\
  0       & L^2       & \nu(y,z)     \\
 \Lambda & \nu(y,z)   & \mu(y,z)
  \end{pmatrix}.
$$
in the moving frame $(\partial_x,\partial_y+nz\partial_x+\theta \partial_z,\partial_z)=(\wdt X,\wdt U, \wdt Z)$. The function $L$ is already constant. We have to see now that we can find new coordinates in which $\nu$ is constant and $L$ is  unmodified. 

 The  lifts  $\wdt U$ and $\wdt Z$ are not unique: according to Lemma \ref{relevons} and \ref{cohomology} we can replace  $\wdt U$ by $\wdt U+\beta (\oU)\partial_x$ and $\wdt Z$ by $\wdt Z + \beta(\partial_z)\partial_x$ where $\beta$ is an \emph{exact} $1$-form on $\T^2$ (in order to keep the same fundamental group, see Lemma \ref{cohomology}).
\begin{fact}\label{koko}
For any smooth function  $\nu:\T^2\rightarrow \R$ there exists a function $N:\T^2\rightarrow \R$ 
such that $\dd N(\oU)=\nu + k$ where $k=-\int\nu\dd\lambda$, with $\lambda$ the Lebesgue measure on $\T^2$. 
\end{fact}
{\bf Proof.} We look at $\T^2$ as the quotient of $\R\times S^1$ by the $\Z$-action defined by identification $1.(s,z)=(s+1,R_\theta z)$. By the previous discussion we can assume that $\oU$ lifts as $\partial_s$.

 Let $N_0 : \R\times S^1\rightarrow \R$ be the function defined by:
$$
N_0(s,z)=\int_0^s \nu(t,z)+ k\ dt. %
$$
Let $h:S^1\rightarrow \R$ be the function defined by:
$$h(z)=N_0(1,z)=\int_0^1 \nu(t,z)+ k\ dt.$$

Clearly, $\partial_s.N_0=\nu+ k $ but $N$ is not invariant by the $\Z$-action. However, for any $\psi:S^1\rightarrow \R$ the function $N_\psi$ defined by $N_\psi(s,z)=N_0(s,z)+\psi(z)$ still satisfies: 
$\partial_s.N_\psi=\nu+k$. Moreover, it will be invariant by the $\Z$-action if and only if $N_\psi(1, z)=N_\psi(0,R_\theta z)$ \ie if $\psi(R_\theta z)-\psi(z)=h(z)$. As $\int_{S^1} h(z) dz=0$, Kolmogorov's Theorem implies that such a function $\psi$ exists. \boxfill$\Box$\medskip

Applying Fact \ref{koko} to the function $-\nu/\Lambda$, we find an exact $1$-form $\beta$ on $\T^2$ such that the expression of $g$ in the moving frame $(\partial_x, \wdt U+\beta (\oU)\partial_x,\wdt Z + \beta(\partial_z)\partial_x)$ is 
$$
\begin{pmatrix}
  0       & 0         & \Lambda\\
  0       & L^2       &  k  \\
 \Lambda &  k   & \mu'(y,z)
  \end{pmatrix}.
$$
According to Lemmas \ref{hophophop} and \ref{cohomology}, there exists coordinates on $\R^3$, that we still denote by $(x,y,z)$, such that $(\partial_{x},\partial_{y}+nz\partial_x+\theta\partial_z, \partial_{z})=(\partial_x, \wdt U+\beta (\oU)\partial_x,\wdt Z + \beta(\partial_z)\partial_x)$ and that the fundamental group of $M$ still reads as $\Gamma_n$. Notice that the coordinates induced on $\T^2$ did not change.
There exists a function $f:S^1\rightarrow \R$ such that, when $\beta$ is replaced by the $1$-form  $\beta+\frac{\dd f}{\dd z}(z)\dd z$, the map   $z \mapsto \int_0^1 \mu'(y,z)\dd y$ is constant, hence equal to ${\cal I}:=\int \mu'\dd \lambda$ with $\lambda$ the Lebesgue measure on $\T^2$. 

In order to have the proper normalization of $k$ and $\mu$, we change the lifts of $\partial_y$ and $\partial_z$.
To keep the form of the matrix of $\wdt g$, the only other possibility is to replace $\partial_y$ by  $\partial_y+p\partial_x$ and $\partial_z$ by $\partial_z+q\partial_x$ with $(p,q)\in \Z^2$. It changes $k$ into $k+(p+\theta q)\Lambda$ and $\mu$ into $\mu+2\Lambda q$. Hence there is only one lift such that $k\in [0,\Lambda[$ and $\int \mu \dd \lambda\in [0,2\Lambda[$.

\label{action_de_phis_et_phit} We suppose now that  $n\neq 0$. 
For any $s\in \R$, the map  $[\Phi^s_{U}:(x,y,z)\mapsto(x+snz, y + s,z)]\in\opna{Diff}(\R^3/\Gamma_n)$ is the time $s$ of the flow of $U=\partial_y+nz\partial_x$. The matrix of $\Phi_{U}^s{}^*\wdt g$ in the moving frame is:
$$
\begin{pmatrix}
  0       & 0         & \Lambda\\
  0       & L^2       &  k +\theta ns\Lambda \\
 \Lambda &  k +\theta ns\Lambda  & \mu'(y+s,z)+2ns\Lambda 
  \end{pmatrix}.$$
  For  $s_0=-\frac{1}{2n\Lambda}{\cal I}$ we have $\int_0^1\mu''(y,z) \dd y=0$, where  $\mu''(y,z)=\Phi^{s_0}_{U}{}^*\wdt g(\partial_z,\partial_z)$. 
In order to kill $k$, we take the pull-back of $\wdt g$ by the flow of $\partial_z$ at time $t_0$: $\Phi_Z^{t_0}:(x,y,z)\mapsto  (x,y,z+t_0) $ where $t_0=\frac k{n\Lambda}+\theta s_0$. The wanted coordinates are built.\medskip

{\em Expliciting the announced group action.} Let $(x,y,z)\mapsto(x',u,v)$ be a change of coordinates preserving the properties stated in Proposition \ref{Diophantine}. As the coefficient $L^2$ in the matrix of $g$ must remain a constant, such a change is necessarily linear on $\pi(M)$ (up to translations, the study of which we postpone) {\em i.e.\@} $(u,v)=(ay+bz,cy+dz)$ with $\text{\footnotesize$\begin{pmatrix}a&b\\c&d\end{pmatrix}$}\in \opna{SL}(2,\Z)$. Conversely, we will show that any such coordinate change on $\pi(M)$ is induced by a coordinate change of $\wdt M$ preserving the properties of Proposition \ref{Diophantine}. Moreover, what precedes shows that $x$ is determined, up to translation, once the coordinates $(y,z)$ are fixed on $\pi(M)$, thus this coordinate change is unique. Let us build it. Take $\text{\footnotesize$\begin{pmatrix}a&b\\c&d\end{pmatrix}$}\in\opna{SL}(2,\Z)$. The corrresponding change $(y,z)\mapsto(u,v)$ is entailed by the following map on $\wdt M$, which normalizes $\Gamma_
n$:
$$\label{formule_changement_coordonnees_sl2Z}(x',u,v)=F(x,y,z)=(x+n(bcyz+\frac{1}{2}ac(y^2-y)+\frac{1}{2}bd(z^2-z)),ay+bz,cy+dz).$$
This map $F$ is close but not equal to the change of coordinates we are looking for.
The matrix of Proposition \ref{Diophantine}, in the coordinates $F(x,y,z)$, is that of $F^{-1\ast}g$ in the moving frame $(\partial_{x'},\partial_{u}+nv\partial_x+\theta'\partial_v,\partial_v)$, where $\theta'$ is such that $\overline{\partial_{u}}+\theta'\overline{\partial_v}$ is tangent to $\overline\FF$ {\em i.e.\@} 
$\theta'=\frac{c+d\theta}{a+b\theta}$. This matrix is equal to   
the matrix of $g$ in the moving frame 
$\mathcal {B} :=(F^{-1})_*(\partial_{x'},\partial_{u}+nv\partial_x+\theta'\partial_v,\partial_z)$. 
Now:
$$F^{-1}_*(\partial_{x'},\partial_{u},\partial_v)=(\partial_{x},(-dB_1+cB_2)\partial_x+d\partial_{y}-c\partial_{z},(bB_1-aB_2)\partial_x-b\partial_{y}+a\partial_z),$$ where $B_1(y,z)=nbcz+nacy-\frac 12 nac$ and $B_2(y,z)=nbcy+nbdz-\frac 12nbd$. Therefore: 
$$
\Mat_{\mathcal {B}}(g)=\left(\begin{array}{ccc}0&0&\Lambda\rho\\
0&L^2/\rho^2&k-L^2b/\rho+\frac n2\Lambda(ac+\theta bd)\\
\Lambda\rho&k-L^2b/\rho+\frac n2\Lambda(ac+\theta bd)&\wdt\mu(u,v)\rho^2+L^2b^2-2kb\rho+n\Lambda ab(d-c)\rho\end{array}\right) 
$$
with $\rho=a+b\theta$ and $\wdt\mu(u,v)=\mu(x,y)=\mu(du-bv,-cu+av)$. We let the reader check the calculations.  
Finally, take $f$ a 1-periodic function. One may turn $x'$ into $x'+f(v)$; this adds $2\Lambda\frac{\dd f}{\dd z}$ to $g_{|(u,v)}(\partial_v,\partial_v)$, hence there is a unique $f$ such that $[v\mapsto\int g_{|(u,v)}(\partial_v,\partial_v)du]={\rm cst}$ {\em i.e.}\@ $[v\mapsto\int \wdt\mu(u,v)\dd u]={\rm cst}$. As ${\cal I}:=\int\!\!\int\wdt\mu\dd u\dd v=\int\!\!\int\mu\dd y\dd z$ remains unchanged, $\wdt\mu$ is turned into  $\mu'=\wdt\mu-{\cal I}(v)$ where ${\cal I}(v):=\int\wdt\mu(u,v)\dd u - \cal I$. 

If $n=0$, we get  
$$
\Mat(g)=
\left(\begin{array}{ccc}0&0&\Lambda\rho\\
0&L^2/\rho^2&k-L^2b/\rho\\
\Lambda\rho&k-L^2b/\rho&\mu'(u,v)\rho^2+L^2b^2-2kb\rho\end{array}\right)\ \text{with }\rho=a+b\theta
$$
{\em i.e.}\@ $(\theta,\Lambda,L,k,\mu)$ is canonically associated with the metric $g$, up to the action of SL$_2(\Z)$ given by $\text{\footnotesize$\left(\begin{array}{cc}a&b\\c&d\end{array}\right)$}.(\theta,\Lambda,L,k,\mu)=(\frac{c+d\theta}{a+b\theta},\Lambda\rho,L/|\rho|,k-L^2b/\rho,\mu'\rho^2+L^2b^2-2kb\rho)$, up to sign for $\Lambda$, $k$ and $\mu$ (see ``{\em Changes of sign}'' below) and up to right    composition by a translation for $\mu$. Notice that the new integral $\int\!\!\int g(\partial_v,\partial_v)\dd u\dd v$ is equal to $\rho^2{\cal I}+L^2b^2-2kb\rho$.

If $n\neq0$, remind that $k=0$ and $\int\!\!\int\mu'\dd u\dd v=\int\!\!\int\wdt\mu\dd u\dd v=\int\!\!\int\mu\dd y\dd z=0$. As a consequence, ${\cal J}:=\int\!\!\int g(\partial_v,\partial_v)\dd u\dd v=L^2b^2+n\Lambda ab(d-c)$. The pull back of $g$ by $\Phi_V^{t_0}\circ\Phi_U^{s_0}$, the flows of the fields $U=\partial_u+nv\partial_x$ and $V=\partial_v$ with:
\begin{equation}\label{s0t0}
\left\{\begin{array}{l}
s_0=-\frac1{2n\Lambda\rho}{\cal J}=-\frac12ab(d-c)-\frac{L^2b^2}{2n\Lambda\rho}\\
t_0=\frac1{n\Lambda\rho}\left(-\frac{L^2b}{\rho}+\frac n2\lambda(ac-\theta bd)\right)+2\theta's\\\phantom{t_0}=-(1+c+d\theta)\frac{L^2b}{n\Lambda\rho^2}+\frac1\rho\left(\frac12(ac+\theta bd)-ab(d-c)(c+d\theta)\right)
\end{array}\right.
\end{equation}
let the last coefficient of the second column of $\Mat(g)$ vanish and adds some constant to its bottom right coefficient so that $\int\!\!\int g(\partial_v,\partial_v)\dd u\dd v$ vanish. As in fact $\int\!\!\int\mu'\dd u\dd v=0$, we get:
$$
\Mat(g)=
\left(\begin{array}{ccc}0&0&\Lambda\rho\\
0&L^2/\rho^2&0\\
\Lambda\rho&0&\mu'(u+s_0,v+t_0)\rho^2\end{array}\right)\ \text{with }\rho=a+b\theta
$$
{\em i.e.}\@ $(\Lambda,L,\mu)$ is canonically associated with the metric $g$, up to the action of SL$_2(\Z)$ given by $\text{\footnotesize$\left(\begin{array}{cc}a&b\\c&d\end{array}\right)$}.(\Lambda,L,k,\mu)=(\Lambda\rho,L/|\rho|,\mu'(u+s_0,v+t_0)\rho^2)$, (see $(s_0,t_0)$ above), up to sign (see just below) and up to right  composition by a translation $(u,v)\mapsto(u+\frac pn,v+\frac qn)$ where $(p,q)\in\Z^2$ for $\mu$. Contrarily to the case $n=0$, those are the only translations in the normalizer of $\Gamma_n$.\medskip

\noindent{\em Changes of sign.} In order to explicit the full action of $\opna{GL}_2(\Z)$ or $\Z/2\Z\times\opna{GL}_2(\Z)$, we look now at the change of $(\theta,\Lambda,L,k,\mu)$ when only the signs of some of the coordinates $x$, $y$, $z$ are changed. If $n=0$, for any $(\varepsilon_1,\varepsilon_2,\varepsilon_3)\in\{\pm1\}^3$, $(x,y,z)\mapsto(\varepsilon_1x,\varepsilon_2y,\varepsilon_3z)$ is in the normalizer of $\Gamma_0$. Its action is: $(\theta,\Lambda,L,k,\mu)\mapsto(\varepsilon_2\varepsilon_3\theta,\varepsilon_1\varepsilon_2\Lambda,L,\varepsilon_2\varepsilon_3k,\mu(\varepsilon_2y,\varepsilon_3z))$. If $n\neq0$, only $(x,y,z)\mapsto(\varepsilon_2\varepsilon_3x,\varepsilon_2y,\varepsilon_3z)$ is in the normalizer of $\Gamma_n$. Its action is: $(\theta,\Lambda,L,\mu)\mapsto(\varepsilon_2\varepsilon_3\theta,\varepsilon_3\Lambda,L,\mu(\varepsilon_2y,\varepsilon_3z))$. \boxfill$\Box$ 

\begin{rem}
The second half of Kolmogorov's Theorem tells us that in general, it is not true that we can find adapted coordinates on $\wdt M$ such that $L$ is constant. It says also that in general we cannot suppose that $\nu$ is constant without changing the coordinates on $\T^2$.  However, it does not say anything about the possibility to find $\nu$ constant thanks to a change of coordinates \emph{and} the choice of appropriate lifts (but such a modification would change the value of $L$).
\end{rem}

\subsection{When the slope is rational.} 

When the leaves of $\overline{\FF}$ are closed, we can also describe  pairwise non isometric models. For the rest of this paragraph $(M,g)$ is a Lorentzian $3$-manifold admitting a parallel lightlike vector field with closed orbits such that the leaves of $\FF$ are also closed.

We start by introducing two new invariants of $(M,g)$ (or at least by giving geometrical interpretation to quantities that already appeared). On $\Mbar=\pi(M)$, each leaf $\overline{\FF}_0$ of $\overline{\FF}$ is a circle, and has a length $L(\FF_0)$ given by $\gbar$. This family of lengths is an invariant of $(M,g)$. Besides, the form $\zeta$ \ie $X^\flat$ pushed on $M/\FF\simeq S^1$ has a total volume $\Lambda :=\int_{M/\FF}\zeta$, which is also an invariant of $(M,g)$. Denoting $\frac1\Lambda\zeta$ by $\dd z$ enables to introduce a ``coordinate'' function $z$ with value in $\R/\Z$, defined, up to translation, on $M/\FF$ as well as directly on $M$. We denote by $F_z$ the leaf of $\FF$ above $z$. The length function $L$ may also be viewed as a function of $z$, setting: $L(z):=L(F_z)$. 

\begin{prop}\label{releve}  There exists $n\in \N$ and a diffeomorphism between $\widetilde M$ and $\R^3$ such that the fundamental group of $M$ is equal to $\Gamma_{n}$ and 
$\wdt g$, the lift of $g$ to $\widetilde M$, reads
$$
\begin{pmatrix}
  0       & 0         & \Lambda\\
  0       & L^2(z)  &  k\\
  \Lambda & k         & \mu(y,z)
  \end{pmatrix}
$$
in the moving frame $(\partial_x,\partial_y+nz\partial_x,\partial_z)$ or equivalently that 
$$ \wdt g=2\Lambda \dd x\dd z+L^2(z)\dd y^2 +2(k-n\Lambda z) \dd y\dd z+\mu(y,z)\dd z^2,
$$
where $\Lambda >0$, $k\in [0,\Lambda[$, $L$ is a $1$-periodic function 
and $ \mu$ is a $(1,1)$-biperiodic function 
such that  $z\mapsto \int_0^1\mu(s,z) ds$ is a constant that belongs to $[0,2\Lambda[$. 
Moreover, if $n\neq 0$ we can assume that $k=0$ and for any $z\in \R,\ \int_0^1 \mu(y,z)\dd y=0$.

Modulo the action of $\Gamma_n$, such coordinates are unique up to the action of:\medskip

-- $((\Z/2\Z)^2\ltimes\Z)\ltimes\R^3/\Z^3$ in case $n=0$,\medskip

-- $(\Z/2\Z\ltimes\Z)\ltimes(\R/\Z\times(\Z/n\Z)^2)$ in case $n\neq0$,\medskip

\noindent where the normal factor acts by translation $(x,y,z)\mapsto(x+x_0,y+y_0,z+z_0)$ with $(x_0,y_0,z_0)\in\R^3/\Gamma_0=\R^3/\Z^3$ if $n=0$ and $(x_0,y_0,z_0)\in\R/\Z\times(n\Z/\Z)^2$ if $n\neq0$. The action of $\Z$ is given by $\toto.(x,y,z)\mapsto(x+\psi_\toto(y,z),y+ \toto H(z) ,z)$  where $H$ is the lift of a diffeomorphism of $\R/\Z$; see the details on $H$ and $\psi_\toto$ in the proof. In particular, it corresponds to the choice of a homotopy class (that of the orbits of $\partial_z$) of a closed transversal to $\overline \FF$ intersecting every leaf once. The factors $\Z/2\Z$ change  the sign of some coordinates; see details at the end of the proof. Therefore, $(\Lambda,L,k,\mu)$ characterize $g$ up to this action. 

The translations preserve $(\Lambda,k)$ and act by natural composition on the functions $(L;\mu)$. The action of $\Z$ preserves $\Lambda$ and $L$; see in the proof its action on $k$ and $\mu$.

The curvature $r=\frac1{L^2\Lambda^2}g(R(Y,Z)Y,Z)$ introduced in Lemma \ref{local} reads: $r=\frac1{\Lambda^2}\left(\frac1{2L^2} \partial^2_y\mu+\frac{\partial_z^2 L}{L}\right)$.
\end{prop}
\begin{rem}
If $n=0$, it is possible to give a slightly different statement of Proposition \ref{releve} by allowing the fundamental group of $M$ to be any of the $\Gamma_{0,c_1,c_2}$ and imposing $k=0$ and for all $z$, $\int_0^1\mu(s,z)ds=0$.
\end{rem}
{\bf Proof.} According to Lemmas \ref{hophophop} and \ref{cohomology}, in order to prove Proposition \ref{releve} we need to find a moving frame $(X,Y,Z)$ satisfying the hypothesis of  Lemma \ref{hophophop} and such that the expression of $g$ in this frame is given by:
$$\Mat_{(X,Y,Z)}(g)=\text{\small$\left(\begin{array}{ccc}0&0&\Lambda \\0&L^2&k\\\Lambda &k&\mu\end{array}\right)$},
$$
with $k\in \R$, $\Lambda \in\R_+^\ast$, $L:M/\FF\rightarrow\R_+^\ast$ and $\mu:M/\mathcal X\rightarrow\R$. 

By integration, the form $\dd z$ defines, up to an additive constant, a function $z:\pi(M)\rightarrow \R/\Z$, whose levels  are the leaves of $\overline{\FF}$. Besides, let $\oY$ be a vector field on $\pi(M)$, tangent to $\overline{\FF}$ and with $\gbar(\oY,\oY)=L^2(z)$ along each $\FF_z$. Actually, this property defines $\oY$ up to its sign. By construction of $\oY$, its flow $\Phi_{\oY}$ is periodic of period $1$. Therefore, taking $\gamma$ any simple loop on $\pi(M)$, intersecting transversely each leaf of $\overline{\FF}$ once, there exists a (unique) function $y:\pi(M)\rightarrow \R/\Z$ such that $y^{-1}(0)=\gamma$ and $\dd y(Y)=1$. Finally, $(y,z)$ is a diffeomorphism $\pi(M)\rightarrow(\R/\Z)^2$ and so $\int_{\pi(M)}\dd y\wedge\dd z=1$. Let $(\oY,\oZ)$ be the basis field dual to $(\dd y,\dd z)$. As $D_ZZ$ and $D_YZ$ are both orthogonal to $X$ for any lift $Z$ of $\oZ$, we have
$g(X,Z)=\Lambda$. Now we modify the fields $Y$ and $Z$ to achieve the wanted form for $\Mat(g)$.

{\mathversion{bold}\bf (i) We get $(Y',Z')$ such that $Y'.g(Y',Z')=0$.\mathversion{normal}} By Lemma \ref{relevons} with $(\Vbar,\Wbar):=(\oY,\oZ)$, we get two fields $(Y,Z)$ on $M$, commuting with $X$, pushed by $\pi$ on $(\oY,\oZ)$, such that $[Y,Z]=nX$; moreover, setting $(Y',Z'):=(Y+\beta(\oY)X,Z+\beta(\oZ)X)$ with $\beta$ any closed 1-form on $\pi(M)$ gives fields with the same properties. Take such a $\beta$, then: $g(Y',Z')=g(Y,Z)+\beta(\oY)\Lambda$. 

We will find a $\beta$ such that $Y'.(g(Y',Z'))=0$. Set $\nu:=g(Y,Z)$ and $N:z\mapsto\int_{\overline{\FF}_z}\nu\dd y$. Then define the form $\beta$ by: $\beta(\oY):=\frac{N-\nu}\Lambda$, $\beta(\oZ)=0$ along some arbitrary simple closed transversal $\gamma_1$ to $\overline{\FF}$, and $\frac{\partial}{\partial y}(\beta(\oZ)):=\frac{\partial}{\partial z}(\beta(\oY))$. To check that this definition of $\beta(\oZ)$ is consistent, we must check that for any $z$, $I_z:=\int_{\overline{\FF}_z}\frac{\partial}{\partial y}(\beta(\oZ))\dd y=0.$ It is the case, as:
$$I_z=\int_{\overline{\FF}_z}\frac{\partial}{\partial z}(\beta(\oY))\dd y=\frac{\partial}{\partial z}\int_{\overline{\FF}_z}\frac{N-\nu}\Lambda\dd y=\frac{\partial}{\partial z}(0)=0.$$ With this $\beta$, $g(Y',Z')=N$ does not depend on $y$.

{\mathversion{bold}\bf (ii) We get $(Y'',Z'')$ such that $g(Y'',Z'')={\rm cst}$.\mathversion{normal}} Take $\tau$ some real number and set: %
$$
     Z'' := Z'-\left(\frac {N}{L^2}+\frac{\tau}{L^2}\right)Y'\ \ \text{and: }Y'' := Y'.\\
$$
As $Y.N$ and $Y.L$ are null, $[Z'',Y'']=[Z',Y']=nX$. Besides, $Y''$ and $Z''$ are still lifts of vector fields $\oY''=\oY$ and $\oZ''$ of $\pi(M)$. As $[\oY'',\oZ'']=0$, $\oZ''$ is conjugated  to the suspension field of some rotation $y\mapsto y+\tau'$, with $\tau'\in\R/\Z$. Choosing adequately the class of $\tau$ modulo $\left(\int_0^1 1/L^2(z)dz \right)\Z$ gives $\tau'=0$. The full choice of $\tau$ gives moreover that the (periodic) orbits of $\oZ''$ are in any of the homotopy class homotopic to $\gamma$ in $\Mbar$.

In other words, we set $Z'':=Z'+h(z)Y'$ with $h$ the unique function from $\R/\Z$ to $\R$ such that $g(Z'',Y'')={\rm cst}$ and that $\int h(z)\dd z$ is the integer such that the orbits of $\oZ''$ are homotopic to $\gamma$.

{\mathversion{bold}\bf (iii) We get $(Y''',Z''')$ such that $z\mapsto\int_0^1\mu(s,z) ds={\rm cst}$.\mathversion{normal}} We remark that if we set $Z'''=Z''+f(z)X$ nothing is changed except that $g(Z''',Z''')=g(Z'',Z'')+2\Lambda f$. It is therefore possible to  find $f$ such that $\int_0^1\mu(s,z) ds={\rm cst}$. We relabel $(Y'',Z''')$ as $(Y,Z)$ and  we are done. 

{\mathversion{bold}\bf (iv) We pull back $\R/\Gamma_{n,c_1,c_2}$ on $\R/\Gamma_{n}$ and put $k$ and $\int\mu(s,z) ds$ to their wished values.\mathversion{normal}} Applying Lemma \ref{hophophop}, we see that  $g$ has the desired expression, but on some  $\R^3/\Gamma_{n,c_1,c_2}$. We take the pull-back of this metric by the diffeomorphism from $\R^3/\Gamma_n$ to $\R^3/\Gamma_{n,c_1,c_2}$ induced by the map $F(x,y,z)=(x+c_1y+c_2z,y,z)$. 
In order to have $k\in [0,\Lambda[$ and $\int_0^1\mu(y,z)\dd y\in [0,2\Lambda[$, we take the pull-back of $g$ by the diffeomorphism of  $\R^3/\Gamma_n$ induced by the map $G(x,y,z)=(x+l_1y+l_2z,y,z)$, where $(l_1,l_2)$ are well chosen elements of $\Z^2$.

If $n\neq 0$, we proceed exactly like at the end of the proof of Proposition \ref{Diophantine} in order to kill $k$ and to have $\int_0^1\mu(y,z)\dd y=0$ for any $z$.

Let us prove the uniqueness of the coordinates modulo the announced group action. Up to sign and translation, the coordinate $z$ is unique. The coordinate $y$, whose levels are the leaves of $\FF$, may be turned into $y'=y-H(z)$ with $H$ such that $H(z+1)-H(z)\in\Z$. We prove that, modulo translation, for each $\toto\in\Z$, there is a unique function $H_\toto$ and a unique coordinate change $(x,y,z)\mapsto(x+\psi_\toto(y,z),y-H_\toto(z),z)$ such that $H_\toto(z+1)=H_\toto(z)+\toto$ and that $\Mat(g)$ has the required form in the moving frame $(\partial_x,\partial_y+nz\partial_x,\partial_z)$.

Take $\toto\in\Z$ and $H_\toto$ such that $H_\toto(z+1)=H_\toto(z)+\toto$. As we postpone the question of the translations, we may suppose that $H_\toto(0)=0$. Then the coordinate change defined by:

$$F_\toto(x,y,z)=(x+n(\frac12-z)H_\toto(z)+n\frac{\toto}{2}z^2,y-H_\toto(z),z)$$
is in the normalizer of $\Gamma_n$. After Lemma \ref{cohomology}, once $(y,z)$ and $z$ are fixed, $x$ may only be turned into $x'=x+\eta(y,z)+Ay+Bz$ with $(A,B)\in\Z^2$ and $\eta$ a (1,1)-biperiodic function. (This change must conjugate $\Gamma_n$ to itself, so corresponds in the Lemma to the choice of a closed form $\beta$ whose cohomology class is entire; here $\beta=-\dd \eta-A\dd y-B\dd z$). Forgetting for the moment the integers $A$ and $B$, {\em any} coordinate change fitting the requirements of Proposition \ref{releve} must be of the form:

$$(x',y',z')=F_\toto(x,y,z)=(x+n(\frac12-z)H_\toto(z)+n\frac{\toto}{2}z^2-\eta(y,z),y-H_\toto(z),z)$$
with $\toto$ and $H_\toto$ as above and $\eta$ some function $\T^2\rightarrow\R$.

To determine which of those changes are actually authorized, we must compute the matrix $G'$ of $(F_\toto^{-1})^\ast g$ in the moving frame $(\partial_{x'},\partial_{y'}+nz\partial_{x'},\partial_{z'})$. We set $h_\toto:=\frac{\dd H_\toto}{\dd z}$ (this function $h_\toto$ plays the same role as the function $h$ introduced above in point {\bf (ii)}). In the frame field $F_{\toto\ast}(\partial_{x},\partial_{y}+nz\partial_{x},\partial_{z})$:

$$\Mat(\partial_{x'},\partial_{y'}+nz\partial_{x'},\partial_{z'})=\text{\small$\left(\begin{array}{ccc}1&\frac{\partial\eta}{\partial y}&n(H_\toto(z)-\frac12h_\toto(z)-\toto z)+\frac{\partial\eta}{\partial z}\\0&1&h(z)\\
0&0&1\end{array}\right)$},$$
hence we get:
$$G'=\text{\small$\left(\begin{array}{ccc}0&0&\Lambda\\
0&L^2&k+h_\toto L^2+\Lambda\frac{\partial\eta}{\partial y}\\
\Lambda&k+h_\toto L^2+\Lambda\frac{\partial\eta}{\partial y}&\mu_1(y',z')+h_\toto^2L^2+2n\Lambda(H_\toto-\toto z-\frac12h_\toto)+2h_\toto k+2\Lambda\frac{\partial\eta}{\partial z}\end{array}\right)$}$$
where $\mu_1(y',z')=\mu(y'+H_\toto(z'),z')$. Hence, $k+h_\toto(z)L^2(z)+\Lambda\frac{\partial\eta}{\partial y}(y,z)$ must be constant {\em i.e.}:\medskip

-- $\frac{\partial^2\eta}{\partial y^2}$ must vanish. As $\eta(y+1,z)=\eta(y,z)$, this implies: $\frac{\partial\eta}{\partial y}=0$. So $h_\toto$ and $H_\toto$ are determined.\medskip

-- $k+h_\toto(z)L^2(z)$ must be a constant {\em i.e.\@} $h_\toto(z)=\frac{\rm Cst}{L^2(z)}$. As moreover $\int_0^1h_\toto(z)\dd z=H_\toto(1)-H_\toto(0)=\toto$, necessarily, {\mathversion{bold}\bf setting ${\cal L}:=1/\int_0^11/L^2(z)\dd z$\mathversion{normal}}, $h_\toto(z)=\frac{\toto{\cal L}}{L^2(z)}$. Thus $H_\toto=\toto H$ with $H(z)=\int_0^z\frac{{\cal L}}{L^2(z')}\dd z'$.\medskip

\noindent Finally, $z\mapsto\int_0^1g(\partial_{z'},\partial_{z'})\dd y$ must be constant {\em i.e.\@} the function $S$ defined by:
\begin{equation}\label{condition_eta}
S(z):=h_\toto^2L^2+2n\Lambda(\toto H-\toto z-\frac12h_\toto)+2h_\toto k+2\Lambda\frac{\partial\eta}{\partial z}\ \text{must be a constant.}
\end{equation}
As $h_\toto$ and $H$ are determined, there is exactly one 1-periodic function $\eta$ achieving this. As  $\int_0^1\frac{\partial\eta}{\partial z}\dd z=0$, the constant value of $S$ is then $\mathcal {J} := \int_0^1 S(z)-2\Lambda\frac{\partial\eta}{\partial z} \dd z$.

Finally, turning $x'$ into $x'+Ay+Bz$ with the adequate integers $(A,B)$, we get that $g(\partial_{y'}+nz\partial_{x'},\partial_{z'})\in\left[0,\Lambda\right[$ and $g(\partial_{z'},\partial_{z'})\in\left[0,2\Lambda\right[$.

If $n=0$, ${\cal J}=\toto(2k+{\toto}{\cal L})$. For any real number $\lambda$, denoting by $\{\,\cdot\,\}_\lambda:\xi\mapsto\lambda\left\{\frac\xi\lambda\right\}$ the ``fractional part modulo $\lambda$'', we get that $\Mat(F_\toto^{-1\ast}g)$ reads, in the frame field $(\partial_{x'},\partial_{y'},\partial_{z'})$:
\begin{equation}\label{matricen0Fferme}
\text{\small
$\left(\begin{array}{ccc}
0&0&\Lambda\\
0&L^2(z')&\left\{k+{\toto{\cal L}}\right\}_\Lambda\\
\Lambda&\left\{k+{\toto{\cal L}}\right\}_\Lambda&\left\{\mu'(y',z')+\toto(2k+{\toto}{\cal L})\right\}_{2\Lambda}
\end{array}\right)$}\ \text{with: }\mu'(y',z')=\mu(y'+H_\toto(z'),z').
\end{equation}
If $n\neq0$, ${\cal J}={\toto^2}{\cal L}+2n\Lambda(-\toto+\int_0^1H_\toto(z)\dd z)$. Turning $(x',y',z')$ into $\Phi_{\partial_{y'}+nz'\partial_{x'}}^{s \ast}\circ\Phi_{\partial_{z'}}^{t\ast}(x',y',z')$ with:
\begin{equation}\label{s_et_t_pour_isometrie}
s\equiv-\frac1{2n\Lambda}{\cal J}=-\frac{\toto^2{\cal L}}{2n\Lambda}+\toto-\int_0^1H_\toto(z)\dd z\ [\frac1n\Z]\quad \text{and: }t\equiv\frac{\toto{\cal L}}{n \Lambda}\ [\frac1n\Z],
\end{equation}
(see the introduction of $s_0$ and $t_0$ in \eqref{s0t0} p.\@ \pageref{s0t0}) we get in the frame field $(\partial_{x'},\partial_{y'}+nz'\partial_{x'},\partial_{z'})$:
$$\Mat(F_\toto^{-1\ast}g)=\text{\small$
\left(\begin{array}{ccc}
0&0&\Lambda\\
0&L^2(z')&0\\
\Lambda&0&\mu'(y',z')
\end{array}\right)$}
\ \ \text{with }\mu'=\mu(y'+H_\toto(z')-s,z'-t).
$$
This final form is obtained only for the given values of $s$ and $t$; the group of admissible translations announced in Proposition \ref{releve} follows.

Eventually, a change of sign is possible. If $n=0$, those in the normalizer of $\Gamma_0$ and preserving $\Lambda>0$ are the: $(\varepsilon_1x, \varepsilon_2y,\varepsilon_1z)$ with $(\varepsilon_1,\varepsilon_2)\in\{\pm1\}^2$. If $n\neq0$, the only one is $(-x, -y,-z)$.

The formula expressing the curvature $r$ 
is a straightforward calculation.\hfill\boxfill$\Box$\medskip

\begin{rem}\label{lairderien}
In case $n=0$ (or if $n\neq0$ and if we fix $z$ in the coordinate change, authorizing $k$ and $z\mapsto \int\mu(y,z)\dd y$ to become non null), it is interesting to note that the Christoffel symbols of the connection are not preserved during the coordinate changes given in Proposition \ref{releve}, {\em except if ${\cal L}\in\Lambda\Q$ and $\toto\in\frac\Lambda{\cal L}\Z$}.
\end{rem}

To conclude this paragraph, we describe the flat affine structure of the leaves of $\FF$ ---~induced by the connection $D$. It turns out that it is encoded in the function $L$.

\begin{prop}\label{holonomy} Let $g$ be the Lorentzian metric on $\R^3/\Gamma_n$ reading in the moving frame  $(\partial_x,\partial_y+nz\partial_x,\partial_z)$:
$$
\begin{pmatrix}
  0       & 0         & \Lambda\\
  0       & L^2(z)  &  k\\
  \Lambda & k         & \mu(y,z)
  \end{pmatrix}.
$$
We denote by $F_z$ the leaf of $\FF$ given by the projection  of $\R^2\times\{z\}$. The leaf $F_z$ is a flat totally geodesic submanifold, hence a flat torus. Its structure, as a flat manifold, is given by a developing map in $\R^2$. The following $\mathcal D$ is a developing map for $(F_z,D)$:
$$\mathcal{D}:(x,y)\mapsto(d_1,d_2)=\left(x-\frac{L(z)\partial_zL(z)}{2\Lambda}y^2,y\right)$$
and its associated holonomy representation is given by:
$$\textstyle\left\{\begin{array}{l}\rho(\gamma_1)=\left[(d_1,d_2)\mapsto(d_1+1,d_2)\right]\\\rho(\gamma_2)=\left[(d_1,d_2)\mapsto\left(d_1-\frac{L(z)\partial_zL(z)}{\Lambda}d_2-\frac{L(z)\partial_zL(z)}{2\Lambda},d_2+1\right)\right].\end{array}\right.$$
where $\gamma_1$ is the map $(x,y)\mapsto (x+1,y)$ and $\gamma_2$ the map $(x,y)\mapsto (x,y+1)$ acting on $\R^2\simeq \widetilde{F_z}$.
\end{prop}
\begin{rem}In fact, as $X$ is parallel and $\partial_y$ is parallel modulo $X$, the linear part of the holonomy group $H$ has to be of the form {\footnotesize$\left\{\left(\begin{array}{cc}1&p\alpha\\0&1\end{array}\right);p\in\Z\right\}$}, thus is characterized by a single real number, $\alpha$. Now this linear part is the only thing that makes sense. Indeed, taking $(\gamma_1,\gamma_2')$ as new basis of $\pi_1(F_z)$, with $\gamma_2'=\gamma_1-\alpha\gamma_1$, we get: $\rho(\gamma'_2):(d_1,d_2)\mapsto(d_1-\frac{L(z)\partial_zL(z)}{\Lambda}d_2,d_2+1)$. Therefore, the structure of $F_z$ is given by the linear part of its holonomy group, itself generated by {\footnotesize$\left(\begin{array}{cc}1&\alpha\\0&1\end{array}\right)$}, written in the basis ${\cal D}(\partial_x,\partial_y)$, and where $\alpha=-\frac1{2\Lambda}\partial_zL^2(z)$.
\end{rem}

\noindent{\bf Proof.} Consider $\mathcal D$ as an identification $\widetilde F_z\simeq \mathcal D(\widetilde F_z)$, then $\partial_{d_1}=\partial_x$ and $\partial_{d_2}=\partial_{y}+\frac{L(z)\partial_zL(z)}{\Lambda}\partial_x$. As $D_{\partial_y}\partial_y=-\frac{L\partial_z.L}{\Lambda}\partial_x$, they are parallel fields of $\widetilde F_z$, hence $\mathcal D$ is a developing map. The rest follows.\boxfill$\Box$

\section{The transformations}\label{transfo}
We want now to describe the manifolds $(M,g)$ such that $\indaf>2$. We do not study here the case $\indaf=2$:  if $(M,g)$ is indecomposable and $[\opna{Aff}(g):\opna{Isom}(g)]>1$ then in fact $[\opna{Aff}(g):\opna{Isom}(g)]=\infty$ by Remark \ref{lemme_phi_preserve_X} just below, and if $(M,g)$ is decomposable, we let this case deliberately apart, see Remark \ref{pasdindicedeux}. So let us take $(M,g)$  indecomposable and such that $\indaf>1$. It follows from sections \ref{reducible} and \ref{notcompact} that $M$ is endowed with a parallel lightlike vector field with closed orbits. Thus, we are going to use the results of sections \ref{closed} and \ref{special} 
 to describe the quotient  $\opna{Aff}(g)/\opna{Isom}(g)$.  

Now $X$ denotes a lightlike parallel vector field with closed orbits and $\varphi$ an element of $\opna{Aff}(g)$.

\begin{rem}\label{lemme_phi_preserve_X} The scalar $\lambda$ such that $\varphi_\ast X=\lambda X$, and $\varphi^\ast X^\flat=\frac1\lambda X^\flat$, given by Lemma \ref{lem_phi_preserve} is $\pm1$, and $\varphi$ preserves the curvature function $r$ defined in Lemma \ref{local}. Indeed, for the first fact, the orbits of $X$ are periodic, and for the second one, Lemma \ref{lem_phi_preserve} ensures that $r\circ\varphi=\lambda^2 r$. It follows also that $\forall n, \varphi^n{}^*g= g+nCX^\flat\otimes X^\flat$, hence that no power of $\varphi$ is an isometry unless $C=0$ {\em i.e.}\@ $\varphi$ is. 
\end{rem}

\begin{lem}\label{dynamic}
If the transformation $\overline \varphi$ is not equicontinuous then the curvature of $\overline D$ is constant along the leaves of $\overline{\FF}$.
\end{lem}
\noindent{\bf Proof.} 
By Remark \ref{lemme_phi_preserve_X}, $\varphi$ preserves $r$. We choose coordinates on $\T^2$ such that $\partial_v$ is tangent to $\overline{\FF}$ and $\zeta(\partial_w)=1$. Let $L$ be the positive function defined by $L^2(v,w)=\gbar_{|(u,v)}(\partial_v,\partial_v)$. For any  $(v,w)\in \Mbar$,
$$d\overline\varphi{}^n(v,w) = \begin{pmatrix}
a_n(v,w) & b_n(v,w)\\
0 &1                        
                        \end{pmatrix}
$$
with $a_n=\frac{L}{L\circ \overline\varphi}$, which is bounded as $M$ is compact.  If $\overline \varphi$ is not equicontinuous then $d\varphi_n(v,w)$, and therefore  $b_n$, tends uniformly to infinity (see \cite{Zeghib}).
We take a subsequence $(\overline\varphi{}^{n_k})$ such that $(\overline\varphi{}^{n_k}(v,w))$  is converging to $(v_\infty,w_\infty)$.

Let $(V_{n_k})$ be the sequence of vectors in the tangent space at the point $(v,w)$ defined by $(\frac{1}{a_{n_k}(v,w)}\partial_v-\frac{1}{b_{n_k}(v,w)}\partial_{w})$. It  is converging to a non zero vector $V_\infty$ tangent to $\overline{\FF}$ and the sequence $(D\overline\varphi^{n_k}(v,w).V_{n_k})$ is converging to $0$. It follows that:
$$r(\overline\varphi^{n_k}(\opna{exp}(V_{n_k})))- r(\overline\varphi^{n_k}(v,w)) = r(\opna{exp}(d\overline\varphi^{n_k}(v,w).V_{n_k}))- r(\overline\varphi^{n_k}(v,w)) \longrightarrow 0.$$
But 
$$r(\overline\varphi^{n_k}(\opna{exp}(V_{n_k})))=r(\opna{exp}(V_{n_k}))\longrightarrow r(\opna{exp}(V_\infty))$$
and 
$r(\overline\varphi^{n_k}(v,w))=r(v,w)$
therefore $ r(v,w)=r(\opna{exp}(V_\infty))$ and $r$ is constant along $\overline{\FF}$. \boxfill$\Box$

\begin{prop}\label{prop_phibar=id}
Let  $(M,g)$ be a compact $3$-manifold endowed with a parallel lightlike vector field $X$ and let $\FF$ be the foliation tangent to $X^\perp$. If $\mathcal X$ and $\FF$ have closed leaves then the diffeormophism $(x,y,z)\mapsto(x+z,y,z)$ in the coordinates of $\widetilde M$ given by Proposition \ref{releve}, defines a diffeomorphism $\varphi\in\opna{Aff}(M,g)\smallsetminus\opna{Isom}(M,g)$. In particular, $\indaf=\infty$. 
\end{prop}

\noindent{\bf Proof.} According to Proposition \ref{releve} we can suppose that $\wdt g=2\Lambda \dd x\dd z+L^2(z)\dd y^2+2(k -nz) \dd y\dd z+\mu(y,z)\dd z^2$ and the fundamental group of $M$ is  $\Gamma_{n,c_1,c_2}$. Let $\wdt \varphi :\R^3\rightarrow \R^3$ be the application defined by $\wdt \varphi (x,y,z)=(x+z,y,z)$. This application is in the normalizer of $\Gamma_{n,c_1,c_2}$ and satisfies $\wdt\varphi{}^*\wdt g=2\Lambda \dd x\dd z+L^2(z)\dd y^2 +2(k-nz) \dd y\dd z+(\mu(y,z)+\Lambda)\dd z^2=\wdt g+\frac{1}{\Lambda}X^\flat\otimes X^\flat$, which is also $D$-parallel. Hence, it induces an application $\varphi: M\rightarrow M$ preserving $D$ and not $g$.\boxfill$\Box$

\begin{rem} In the situation of Proposition \ref{prop_phibar=id}, $\opna{Aff}(g)$ contains other types of non isometric elements. You may see them all appear at the end of the proof of Theorem \ref{theo_equivariant}.  
\end{rem}
When the leaves of $\FF$ are not closed the situation is more complicated. We prove the following.

\begin{prop}\label{transfo irra}

Let $(M,g)$ be a compact, orientable, $3$-dimensional Lorentzian manifold admitting a lightlike parallel vector field $X$ such that the leaves of $\mathcal X$ are closed and  those of $\FF$ are not. 
There exists a diffeomorphism preserving $D$ and not $g$ if and only if either $g$ is flat or $(M,g)$ is isometric to $(\R^3,\wdt g)/\Gamma_{n}$ with $n\neq 0$ and $\wdt g$ reads
$$
\begin{pmatrix}
  0       & 0         & \Lambda\\
  0       & L^2       & 0 \\
  \Lambda & 0        & \mu(y)
  \end{pmatrix}
$$
in the moving frame $(\partial_x,\partial_y+nz\partial_x+\theta \partial_z,\partial_z)$ or equivalently
$$
\wdt g= 2\Lambda(\dd x-nz\dd y)(\dd z-\theta\dd y)+L^2\dd y^2+\mu(\dd z-\theta\dd y)^2,
$$
where $\Lambda$ and $L$ are in $\R^*_+$, $\theta\in (\R\smallsetminus \Q) \cap\, ]0,1[
$ and $\mu$ is a $1$-periodic function such that $\int_0^1\mu(y)\dd y=0$. 
Then, the coordinates $(x,y,z)$ of $\widetilde M$ given above are unique, up to translation $(x,y,z)\mapsto(x+x_0,y+\frac pn,z+\frac qn)$ with $x_0\in\R$ and $(p,q)\in\Z^2$. Hence, $(n,\theta,\Lambda,L,\mu)$, up to a change of $\mu$ into $\mu(\,\cdot\,+\frac pn)$, characterize the isometry class of $g$.

In case $M$ is a flat parabolic torus, it has lightlike fibre and $\opna{Aff}(g)/\opna{Isom}(g)$ is spanned by the class of  a flow $(\Phi^t_Y)_{t\in\R}$ where $Y$ is a certain spacelike vector field orthogonal to this fibre. Else it is spanned by the class of the diffeormorphism $\varphi_0$ (well)-defined by its action $\widetilde\varphi_0:(x,y,z)\mapsto(x+z,y,z+\frac\theta n)$ on $\widetilde M$.
\end{prop}
{\bf Proof.} 
First we assume that $g$ is flat. If $M$ is a $3$-torus, it is complete according to Carri\`ere's Theorem \cite{carriere2} and therefore it is the quotient of Minkowski space by translations. It follows that (a conjugate of) the whole group $\rm{GL}(3,\Z)$ acts on $M$ preserving the connection and it is not difficult to find an element that is not an isometry (even if it can be tricky to determine $\rm{GL}(3,\Z)\cap \rm{Isom}(g)$).
We assume now that $M$ is not a torus, then it is a parabolic circle bundle over $\T^2$ . It is well-known (see \cite{Goldman})  that flat Lorentzian metrics on parabolic circle bundles are quotients by a lattice of the Heisenberg group $H$ endowed with a left invariant metric. Moreover, left invariant metrics on $H$ are flat if and only if the center of the Lie algebra of $H$ is lightlike. As we have already seen, there exists  diffeomorphisms between the Heisenberg group and $\R^3$  sending  a moving frame of left  invariant vector fields of $H$ on ($\partial_x, \partial_y+nz\partial_x, \partial_z)$.
Otherwise said there exists coordinates on $\R^3$ and $(a,b,c,d,e)\in \R^5$ such that $\wdt g$, the lift of $g$ to the universal cover of $M$ reads $$\wdt g=a(\dd x-nz \dd y)\dd y +  b(\dd x-nz \dd y)\dd z+ c \dd y^2 + d \dd y\dd z + e \dd z^2$$ 
and that the fundamental group is $\Gamma_{n}$ with $n\neq 0$ (up to isometry we could also assume that $d=e=0$). The vector field $\partial_y+nz\partial_x -\frac {a}{b} \partial_z $ induces a vector field $Y$ on $M$ that preserves $D$ and not $g$. More precisely, we have $\Phi_Y^t{}^*g=g+\frac{2nt}{b}X^\flat\otimes X^\flat$, where $ \Phi_Y^t$ is the flow of $Y$. %
It follows from Lemma \ref{lemme_phi_preserve_X} that any diffeomorphism preserving $D$ is the composition of an isometry of $g$ and  a $\Phi_Y^{t}$.\hip

We assume now that $g$ is not flat. Let $\varphi$ be an element of $\rm{Aff}(g)$. 
In this case  $\overline \varphi$, the application induced on $\Mbar$ by $\varphi$, has to be equicontinuous.
 Indeed, else the curvature of $\overline D$ would be  constant by Lemma \ref{dynamic}  and therefore, by Proposition \ref{GaussBonnet}, $g$ would be flat.

The fact that the  application $\overline\varphi$ is equicontinuous and preserves a connection implies that it preserves some Riemannian metric. It means that  either $\overline \varphi$ has  finite order or the set $\{\overline\varphi^n, n\in \Z\}$ is dense in some 1- or 2-dimensional torus of $\opna{Aff}(\overline D)$. We show that if $\varphi$ is not an isometry, then $\overline\varphi$ is of infinite order and $n\neq0$; by the way, we build the announced coordinates $(x,y,z)$ of $\widetilde M$. Finally we determine the expression of $\varphi$.

Let us first suppose that there exists $k>0$ such that $\overline \varphi^k=\rm{Id}$.
According to Proposition \ref{enforme}, we can identify $M$ with $\R^3/\Gamma_n$ and $\wdt g$, the lift of $g$ to $\wdt M$, reads:
$$
\wdt g= 2\Lambda(\dd x-nz\dd y)(\dd z-\theta\dd y)+L^2\dd y^2+2\nu\dd y(\dd z-\theta\dd y)+\mu(\dd z-\theta\dd y)^2,
$$
where $L$, $\nu$ and $\mu$ are functions.
We choose a lift of $\varphi^k$ to $\R^3$, denoted by $\wdt\varphi{}^k$, that  preserves the coordinates $y$ and $z$.  Hence, there exists two functions $\alpha$ and $\beta$ such that $$(\wdt\varphi{}^k)^*\dd x=\dd x+\alpha\dd y +\beta \dd z.$$
Thus, as the functions $L$, $\nu$ and $\mu$ do not depend on $x$ and as $\varphi^*dy=dy$ and $\varphi^*dz=dz$:
$$(\wdt\varphi{}^k)^*\wdt g=\wdt g+2\Lambda\left( -\theta\alpha\dd y^2 + (\alpha -\theta\beta)\dd y\dd z + \beta \dd z^2\right).$$
Moreover, we know that  $\wdt X{}^\flat=-\Lambda\theta\dd y+\Lambda \dd z$ and that there exists $C\in \R$ such that $(\wdt\varphi{}^k)^*\wdt g=\wdt g +C\wdt X{}^\flat\otimes \wdt X{}^\flat= \wdt g + C\Lambda^2(\theta^2 \dd y^2 -2\theta \dd y \dd z +\dd z^2)$. It follows that $\alpha$ and $\beta$ are constant and that $\alpha=-\beta \theta$ (it means in particular that $\varphi^k$ preserves a flat metric on $M$). Without loss of generality, we can suppose that $\wdt \varphi $ fixes $0$. It means that 
$$\wdt\varphi{}^k(x,y,z)=(x-\beta\theta y+\beta z,y,z).$$

The fact that $\widetilde {\varphi^k}$ lies in the normalizer of $\Gamma_{n}$ implies that $\beta\in \Z$ and $\beta\theta\in \Z$. As $\theta\not\in \Q$ it follows that $\alpha=\beta=0$ \ie $\varphi^k=\opna{Id}$ and therefore $\varphi$ is an isometry.\hip

Now we suppose that the set $\{\overline\varphi{}^n, n\in \Z\}$ is dense in a torus of $\opna{Aff}(\overline D)$. We are not going to use Proposition \ref{enforme} but rather prove a more adapted version. If this torus is at least $2$-dimensional, it would act transitively on $\overline M$ so $\overline M$ would be flat, which is not the case.
 It means that $\overline\varphi$ is the time $t_0$ of a $1$-periodic flow preserving $\overline D$. We denote by $\overline K$ the vector field associated with this flow that satisfies $\zeta(\overline K)>0$. 

Let $\overline Y$ be a parallel vector field, unique up to sign, tangent to $\overline \FF$ (\ie  having a lift on $M$ of constant ``norm'') such that the integral of $\omega$,  the $2$-form dual to the frame $(\overline Y,\overline K)$, on $\mathbb T^2$ is equal to $1$. As $\overline \varphi$ is equicontinuous and volume preserving we see that $\overline\varphi_*\overline Y=\oY$ and therefore $[\overline K,\oY]=0$. Let us choose the  lift $K$ of $\overline K$ to $M$ which is lightlike and  the  lift $Y$ of $\overline Y$ which is orthogonal to $K$.  As $Y$ and $K$ are lifts of $\oY$and $\overline K$, and as $[\oY,\overline K]=0$, there are functions $f_1$, $f_2$ and $f$ such that $[X,Y]=f_1X$, $[X,K]=f_2X$ and $[Y,K]=fX$. Moreover, $X$ is $g$-Killing and $Y$ and $K$ are the only lifts of $\oY$ and $\overline K$ such that  $g(K,K)\equiv0$ and $Y\perp\Span(X,K)$, so $f_1=f_2=0$ and $[X,[Y,K]]\equiv0$, which gives that $X.f=0$.

Computing the curvature of $g$ we find 
$$R(K,Y)K={\rm Cst}\,(Y.f)Y.$$
Thus the function $r$  characterizing the curvature of $\overline D$ defined in Proposition \ref{GaussBonnet}, is equal to $\rm{Cst}\,(Y.f)$.
But clearly $r$ is constant along the integral curves of $\overline K$, therefore we have $K.Y.f=0$ and so $Y.K.f=0$. It means that $K.f$ is constant along the leaves of $\FF$ which are dense. Thus $K.f=0$. Consequently, there exists a function $F$ such that $Y.F=f+ \mrm{Cst}$. Replacing $K$ by $Z=K-FX$ we have $[Y,Z]= \mrm{cst}.X$.  
It follows from the fact that $\int_{\T^2}\omega =1$ that $\mrm{cst}=n$ (see the proof of Lemma \ref{relevons}). 
There exits coordinates $(y,z)$ on $\T^2$ and $\theta\in \R\smallsetminus \Q$ such that $\partial_z=\overline{Z}$  and $\oY=\partial_y+\theta\partial_z$. We may require moreover that $\theta\in\left]0,1\right[$; this fixes $\partial_y$, which else is defined up to addition of $p\overline{Z}$, $p\in\Z$.
Applying Lemma \ref{hophophop}, we obtain an identification between $M$ and $\R^3/\Gamma_{n,c_1,c_2}$ such that:
$$
\wdt g= 2\Lambda(\dd x-nz\dd y)(\dd z-\theta\dd y)+L^2\dd y^2+\mu(y)(\dd z-\theta\dd y)^2.
$$
We notice that this time $L$ is a constant, because $Y$ is parallel modulo $\langle X\rangle$, and $\mu$ depends only on $y$, because $Z.f=0$. As $F$ is only defined up to a constant we can impose $\int_0^1 \mu(y)\dd y=0$. Moreover if $n=0$ then $\varphi$ is an isometry. Indeed, in this case $K$ is Killing and we can suppose $\overline \varphi=\operatorname{Id}$, but we just proved that it means that $\varphi$ is an isometry.

We are left with the case  $n\neq 0$. In order to  work  on $\R^3/\Gamma_n$  instead of $\R^3/\Gamma_{n,c_1,c_2}$, we take the pull-back of $g$ by the diffeomorphism induced by $G(x,y,z)=(x+c_1y+c_2z,y,z)$. The matrix associated  with  $G^*\wdt g$ is 
$$\begin{pmatrix}
0       & 0         & \Lambda\\
0       & L^2       &  \Lambda(c_1+\theta c_2) \\
\Lambda & \Lambda(c_1+\theta c_2)         & \mu(y) + 2\Lambda c_2
\end{pmatrix}.
$$ 
In order to kill the terms $\Lambda(c_1+\theta c_2)$ and $2\Lambda c_2$ that appeared during this operation, we proceed as at the end of the proof of Propositions \ref{Diophantine} and \ref{releve}. With the same kind of arguments as at the end of the proof of Propositon \ref{Diophantine}, it follows from all the construction that the coordinates $(x,y,z)$ are unique up to the announced translations.

Let  $f:\R^3\rightarrow \R^3$ be the map defined by $f(x,y,z)=(x+nt_0y,y,z+t_0)$, where $t_0$ is the real number introduced above. This map is an isometry of $\wdt g$:  this  follows from the fact that the moving frame $(\partial_x,nz \partial_ x  +\partial_y+\theta\partial_z,\partial_z)$ and the coordinate $y$ are invariant by $f$.
Possibly composing  with a translation in the direction of $\partial_x$, we can suppose that the map $f^{-1}\circ \wdt\varphi$ fixes $0$. As it fixes also the coordinates $y$ and $z$, it follows from the previous case  that there exists $\beta\in \R$ such that  $f^{-1}\circ \varphi(x,y,z)=(x - \theta \beta y + \beta z,y,z)$. Hence $\varphi(x,y,z)=(x+(n t_0-\theta \beta)y+\beta z, y, z+t_0)$. 

By a direct computation, we see that $\widetilde \varphi$ is in the normalizer of the fundamental group if and only if $\theta\beta-nt_0\in\Z$ and $\beta\in \Z$, therefore $\theta\beta-nt_0=0$ (check it).  Hence, our map $\widetilde \varphi$ is a power of the map $\widetilde \varphi_0$ defined by  
$$\widetilde \varphi_0(x,y,z)= (x+z,y,z+\theta/n).$$
By a direct computation, we check that $\varphi_0$ is indeed an element of $\opna{Aff}(g)\smallsetminus\opna{Isom}(g)$.
\boxfill$\Box$\medskip

\begin{prop}\label{champs}
Let $(M,g)$ be a $3$-dimensional compact Lorentzian manifold.  Then $\dim(\opna{Aff}(g)/\opna{Isom}(g))\leqslant 1$ with equality if and only if
there exists $n\in \N\smallsetminus\{0\}$ and a diffeomorphism between $\widetilde M$ and $\R^3$ such that the fundamental group of $M$ is equal to $\Gamma_{n}$ and either $\wdt g$ reads:
$$
\wdt g=2\Lambda (\dd x-nz\dd y)\dd z+L^2(z)\dd y^2,
$$
where $L$ is a $1$-periodic function ---~then ${\FF}$ has closed leaves~---, or  $\wdt g$ reads:
$$\wdt g= a(\dd x-nz \dd y)\dd y +  b(\dd x-nz \dd y)\dd z+ c \dd y^2 $$
where $(a,b,c)\in \R^3$,
\ie $(M,g)$ is a flat parabolic torus ---~be ${\FF}$ with closed leaves or not.

\medskip

Otherwise stated, there exists a vector field $K$ whose flow preserves $D$ and not $g$ if and only if $(M,g)$ admits a lightlike parallel  vector field with closed  orbits, $M$ is not a torus and the curvature fonction $r$ defined in Proposition \ref{GaussBonnet} is constant along $\FF$. Then $\opna{Aff}(g)/\opna{Isom}(g)$ is spanned by a  vector field $K=\partial_y+nz\partial_x+\theta \partial_z$ with $\theta=0$ when $g$ is not flat and $\theta=-a/b$ in the flat case. 
\end{prop}
{\bf Proof.} If there exists a vector field $K$ preserving $D$ and not $g$ then  it follows from sections \ref{reducible} and \ref{notcompact} that there exists on $(M,g)$ a parallel lightlike field $X$ with closed  orbits.
We know also that $(M,g)$ is not a flat torus: the Lorentzian flat tori are, affinely, $\R^3/\Z^3$, thus their affine group is $\opna{GL}(3,\Z)$, which is discrete. 

It follows then from Proposition \ref{decomposable} that $(M,g)$ is indecomposable and therefore for any $\varphi\in \opna{Aff}(g)$ we must have $\varphi^*g=g+CX^\flat\otimes X^\flat$ for some number $C$. Moreover,   if $\varphi'$ is another element of $\opna{Aff}(g)$ such that $\varphi'{}^*g=g+CX^\flat\otimes X^\flat$ then $\varphi'\circ \varphi^{-1}$ is an isometry of $g$. Consequently, $\dim(\opna{Aff}(g)/\opna{Isom}(g))\leqslant 1$.

\begin{lem}\label{colin}
The vector field $K$ is not pointwise collinear to $X$.
\end{lem}

\begin{rem}This kind of situation actually happens on the universal cover.
\end{rem}
\noindent{\bf Proof.}
Let us suppose that $K$ is pointwise collinear to $X$. It means that there exists a function $h$ such that  $K=hX$ and $X.h=0$ (otherwise $K$ is not bounded). The vector field $X$ being equicontinuous, it is not difficult to find a moving frame $(X,Y',Z')$ such that:
$$\operatorname{Mat}_{(X,Y',Z')}(g)=\text{\small$\begin{pmatrix}
                                    0&0&1\\0&1&0\\1&0&0
                                    \end{pmatrix}$},\ \ [X,Y']=0\text{ and }[X,Z']=0.  
$$
We have $[K,X]=0$, $[K,Y']=-(Y'.h)X$ and $[K,Z']=-(Z'.h)X$.
 As ${\cal L}_Kg={\rm Cst}.X^\flat\otimes X^\flat$,  we see that $Y'.h=0$ and $Z'.h$ is constant. But as $M$ is compact, it implies that $h$ is constant and therefore $K$ preserves $g$. Impossible.\boxfill$\Box$
\hip
It follows from Remark \ref{lemme_phi_preserve_X} that  $\mathcal L_K X^\flat=0$ and therefore, as $\dd X^\flat=0$, $X^\flat(K)$ is constant.

 By Lemma \ref{basic_pres}, $\Phi_K^t{}^*g (K,K)=g(K,K)+ {\rm Cst.} tg(K,X)^2$. Moreover,  
$\Phi_K^t{}^*g (K,K)=g(K,K)\circ \Phi_K^t$. As $M$ is compact, the function $g(K,K)$ is bounded, thus $g(X,K)=0$, thus $K$ is tangent to $\FF$. 

The vector field $K$ induces a vector field $\overline K$ on $\overline M$. This vector field preserves $\overline D$. By Lemma \ref{colin}, $\overline K$ is non zero on a open dense subset of $\overline K$ (a vector field preserving a connection is determined by its $1$-jet at a point).
It implies that the curvature of $\overline D$ (or more precisely the function $r$ defined in \ref{GaussBonnet}) is constant along the leaves of $\overline {\FF}$.
If these leaves  are dense in $\overline M$ then the curvature of $\overline D$ is constant, and equal to $0$ by the Gauss-Bonnet relation of Proposition \ref{GaussBonnet}. 

If the leaves of $\FF$ are closed, by Proposition \ref{releve} there exists coordinates on $\R^3\simeq \widetilde M$ such that 
$$
\wdt g=2\Lambda \dd x\dd z+L^2(z)\dd y^2 +2(k-\Lambda nz) \dd y\dd z+\mu(y,z)\dd z^2,
$$
Moreover we know by Proposition \ref{releve} that  the curvature $r$ of $\overline D$ is given by $r=\frac1{\Lambda^2}\left(\frac1{2L^2} \partial^2_y\mu+\frac{L''}{L}\right)$, therefore, as $\partial_yr=0$,  $\partial_y\mu=0$. But $\mu$ is such that  $ z\mapsto\int_0^1\mu(s,z)ds$ is constant, therefore $\mu$ is constant.
  As, in a moving frame $(\partial_x,\partial_y+nz\partial_x,\partial_z)$:
$$\Mat(\widetilde g)=\text{\small$\left(\begin{array}{ccc}0&0&\Lambda\\0&L^2&k\\\lambda&k&\mu\end{array}\right)$}\ \text{and:\ }\Mat(X^\flat\otimes X^\flat)=\text{\small$\left(\begin{array}{ccc}0&0&0\\0&0&0\\0&0&\Lambda^2\end{array}\right)$},$$
and as $\Phi_K^{t\ast}g=g+{\rm Cst.}tX^\flat\otimes X^\flat$, we can assume, replacing $g$ by an adequate $\Phi_K^{t\ast}g$ ---~which is by construction isometric to $g$~---, that $\mu=0$. 

Let us prove now  that  $n\neq 0$. Let 
 $g_0$ be the flat metric defined by $\wdt g_0=2\Lambda \dd x\dd z+ \dd y^2 - 2\Lambda nz \dd y\dd z$.
We see that $K$  also preserves the connection of $g_0$ and that $K$ preserves $g_0$ if and only if $K$ preserves $g$. But any affine flow on a flat Lorentzian torus is a Killing flow. 
Therefore $n\neq 0$ and, by Proposition \ref{releve}, we can assume that $k=0$.
\hop

Reciprocally,  a direct computation  proves that if $n\neq 0$ and  there exists coordinates $\R^3\simeq \widetilde M$ such that 
$$
\wdt g=2\Lambda (\dd x-nz\dd y)\dd z+L^2(z)\dd y^2,
$$
then the vector field $\partial_y+nz\partial_x$, which is invariant by $\Gamma_{n}$, preserves $D$ and not $g$. 

Finally, the flat case  with irrational slope of $\overline{\FF}$ in $\overline M$  follows from the discussion at the beginning of the proof of Proposition \ref{transfo irra}.
 \boxfill$\Box$\medskip

To conclude we prove the following result announced in the introduction:
\begin{te}\label{theo_equivariant}
Let $(M,g)$ be a compact orientable and time-orientable Lorentzian $3$-manifold such that  $[\opna{Aff}(g):\opna{Isom}(g)] >2$.
Then for any $\varphi\in\opna{Aff}(g)\smallsetminus \opna{Isom}(g)$, respectively $\varphi\in\opna{Isom}(g)$, there exists a smooth path of metrics $g_t$ on $M$  between $g$ and a flat metric 
 satisfying 
$\varphi\in \opna{Aff}(g_t)\smallsetminus \opna{Isom}(g_t)$, respectively $\varphi\in \opna{Isom}(g_t)$, for any $t$.

However, if $g$ admits an affine non-isometric automorphism  inducing a non equicontinuous transformation on $\overline M$, if $M$ is not a torus and if $g$ is not flat ---~such metrics exixst---, then $\opna{Aff}(g)$ is not included in the affine group of any flat metric.
\end{te}

\begin{rem}The metrics of the second statement of Theorem \ref{theo_equivariant} are those of case (6) with $n\neq0$ of Table \ref{table2}. We see then that they are exacty the non flat metrics $g$ on parabolic tori such that $\afsuris>2$ and $\opna{Isom}(g)$ is non compact.
\end{rem}

{\bf Proof.} We suppose that $g$ is not flat otherwise there is nothing to prove.
As $[\opna{Aff}(g):\opna{Isom}(g)] >2$, there exists on $(M,g)$ a parallel lightlike field with closed  orbits. We can apply Proposition \ref{enforme} in order to have a description of $(M,g)$. 
It follows from Propositions \ref{transfo irra} and \ref{releve} that the function $\nu$ appearing in the description can be supposed  to be null. In most cases, the path of metrics is given by any affine path between the function $L$ and a constant and between the function $\mu$ and a constant. More precisely, $g_t$ is the metrics on $\R^3/\Gamma_n$ whose matrix in the moving frame $(\partial_x,\partial_y+nz\partial_x+\theta \partial_z,\partial_z)$ (the parameter $\theta$ being given by the metric $g$) is 
\begin{eqnarray}
\begin{pmatrix}
0       & 0        &  \Lambda\\
0       & (1-t) + t\, L^2(y,z)   &  k\\
\Lambda & k & t\, \mu(y,z)
 \end{pmatrix}, 
\end{eqnarray}

Let $\varphi$ be an element of $\opna{Aff}(g)\smallsetminus \opna{Isom}(g)$. 
If the leaves  of ${\FF}$  are not closed then by Proposition \ref{transfo irra} and its proof we know that we can choose $(x,y,z)$ such that  $\varphi$  preserves the coordinate $y$ and the expression of the  metrics does not depend on the coordinate $z$. Then it can be checked directly that it implies that $\varphi$  lies in $\opna{Aff}(g_t)\smallsetminus \opna{Isom}(g_t)$ for all $t$.

If the leaves of $\FF$ are closed, let us consider $\varphi\in\opna{Aff}(g)\smallsetminus \opna{Isom}(g)$ and $\overline\varphi$ its action on $\overline M$. If $\overline\varphi=\Id$ then $\varphi$ acts as $(x,y,z)\mapsto(x+z,y,z)$ on $\widetilde M$ in the coordinates of Proposition \ref{prop_phibar=id}, see this proposition and its proof. Then $\varphi\in\opna{Aff}(g_t)\smallsetminus \opna{Isom}(g_t)$ for all $t$. Four other possibilities may occur for $\overline\varphi$; {\bf (i)}-{\bf (iii)} form the case where $\overline\varphi$ is equicontinuous, hence preserves some Riemannian metric. 
\begin{itemize} 
\item[{\bf (i)}]
Either $\overline\varphi$ is periodic. According to Proposition \ref{releve} and its proof, applied with $\toto=0$, $\wdt\varphi=\tau_{p,t}\circ\tau'_{q,s}$ with $\tau'_{q,s}:(x,y,z)\mapsto (x+(q+ns)z,y+s,z)$ where $q\in\Z$ and, 
$s$ arbitrary in $\Q\cap[0,1[$,
 and with $\tau_{p,t}:=(x,y,z)\mapsto(x+py,y,z+t)$ with $p=0$ and $t$ arbitrary in $\Q\cap[0,1[$ if $n=0$, else $p\in\llbracket0,n-1\rrbracket$ and $t=\frac pn$. (In fact, from Proposition \ref{releve} we get the form of $\wdt\varphi$ if $\varphi$ is {\em isometric}, and for them the constraint on $(q+ns,s)$ is stronger: it is the same as that on $(p,t)$ above. We let the reader check that considering affine applications gives the above looser constraint on $q$ and $s$.) 

Finally, with the introduced notation,  for some $q\in\Z$, $\wdt\varphi(x,y,z)=(x+qz,y+s,z+t)$ and $(s,t)\in(\Q\cap[0,1[)^2$ if $n=0$, else $\wdt\varphi(x,y,z)=(x+py+(q+ns)z,y+s ,z+\frac pn)$.

Let $r$ be the smallest integer such that $rt\in\N^\ast$. By Proposition \ref{holonomy}, $L$ is necessarily $\frac1r$-periodic. The fact that $\varphi$ is affine implies also that there is some constant $C$ such that $\mu(\overline\varphi(y,z))=\mu(y,z)+C$; then necessarily $C=0$ as $\overline M$ is compact. Conversely, if $L$ and $\mu$ satisfy those periodicities, the maps $\varphi$ defined above are indeed $g$-affine. Now all those $\varphi$ preserve the canonical flat affine structure of the Heisenberg group if $n\neq0$, or of $\R^3$ if $n=0$ \ie  is affine for the flat Lorentzian metric $g_0$ given above, as well as for each of the $(g_t)_{0\leqslant t\leqslant 1}$.

\item[{\bf (ii)}]
Either  $\{\overline \varphi^l, l\in \Z\}$ is dense in a $1$-dimensional subgroup of $\opna{Diff}(M)$ whose orbits are the leaves of $\overline \FF$.
 This vector field being affine and equicontinuous, it is proportional to $\partial_y$.   The curvature of $\oD$ is  then  constant along the leaves of $\overline{\FF}$ and the metric has the  form  given by Proposition \ref{champs}. In this case, the vector field $Y=\partial_y+nz\partial_x$ is affine (Killing if $n=0$) and tangent to $\overline \FF$ therefore $\varphi=\Phi_Y^s\circ \varphi_0$ with $\overline \varphi_0=\opna{Id}$, and $\varphi\in\opna{Aff}(g_t)\smallsetminus \opna{Isom}(g_t)$ for all $t$.

\item[{\bf (iii)}] Either  $\{\overline \varphi^l, l\in \Z\}$ is dense in a $1$-dimensional subgroup of $\opna{Diff}(M)$ whose orbits form a foliation $\overline {\mathcal G}$ transverse to $\overline \FF$. As we saw in the proof of Proposition \ref{transfo irra}  we can suppose that $\partial_z$ is tangent to $\overline {\mathcal G}$.
It implies that the holonomy of all the tori $F_z$ are equal. Therefore, by Proposition \ref{holonomy}, the function $L$ is constant. The curvature being constant along the orbits of  $\overline \varphi$  we see also that $\partial_z\mu=0$.  Hence the metric does not depend on the coordinate $z$ and $\varphi$ preserves the coordinate $y$. Again, it is enough to conclude that $\varphi\in\opna{Aff}(g_t)\smallsetminus \opna{Isom}(g_t)$ for all $t$. But we can give a better description, indeed it follows from Proposition \ref{releve}, that necessarily $n=0$ and therefore $\partial_z$ is Killing. It implies that $\varphi=\Phi_{\partial_z}^s\circ \varphi_0$ with $\overline\varphi_0=\opna{Id}$. 
\item[{\bf (iv)}] Or $\overline \varphi$ is not equicontinuous. 
According to Lemma \ref{dynamic} and the proof of Proposition \ref{champs} there exists coordinates  $(x,y,z)$ of $\wdt M$ such that $\wdt g$ reads $2\Lambda(\dd x-nz\dd y)\dd z +L^2(z)\dd y^2 +2k_1\dd y\dd z+ k_2 \dd z^2$ with $k_in=0$ for $i=1$ and $2$.
As $\varphi$ preserves $X^\flat$ it acts by translation on the coordinate $z$. Hence, the orbit of any point of $M/\FF$ is finite or dense. If it is dense then all the leaves of $\FF$, as flat affine tori,  have the same holonomy and, by Proposition \ref{holonomy}, $(M,g)$ is flat.

According to Proposition \ref{releve} and its proof, if $\varphi$ is affine, then $\mathcal L/\Lambda$, where $\mathcal L=1/\int_0^1\frac{1}{L^2(s)}\dd s$, must be rational  (see Remark \ref{lairderien}) and $\wdt\varphi=\tau\circ\varphi_{\toto,A,B}$ with:\medskip

-- $\tau$ an affine morphism as described in point {\bf (i)}, except that only its action on the coordinate $z$ has to be periodic, \ie that $s$ may be any real number in $\left[0,1\right[$ instead of a rational one,\medskip

-- $\varphi_{\toto,A,B}$ the affine morphism fixing $0$ defined by:
\begin{equation}\label{eq_philab}
\wdt\varphi_{\toto,A,B}(x,y,z)=(x+n\toto(\frac12-z)H(z)+n\frac{\toto}{2}z^2-\eta(z)+Ay+Bz,y-\toto H(z),z)
\end{equation}

where $(\toto,A,B)\in\Z^3$, where $H(z)={\mathcal L} \int_0^z \frac{1}{L^2(s)}\dd s$, where $\eta$ is defined by Condition \eqref{condition_eta} p.\@ \pageref{condition_eta}, where $\toto\neq0$ (else $\varphi$ would be equicontinuous) and where $A\neq0$ and $A=-\frac{\toto{\cal L}}\Lambda$.\medskip

\noindent {\em Remark.} For a given $\ell$, the $\varphi_{\toto,A,B}$ differ by a translation. Of course, some of the $\varphi_{\toto,A,B}$ {\em e.g.\@} the $(\varphi_{\ell,\ell{\cal L}/\Lambda,0})_{\ell\in\Z}\simeq\Z$ form a group.\medskip

Finally, $\wdt\varphi$ reads $(x,y,z)\mapsto(x-\eta(z)+Ay+Bz,y-\toto H(z)+s,z+ t)$ if $n=0$, else $(x,y,z)\mapsto(x+n\toto(\frac12-z)H(z)+n\frac{\toto}{2}z^2-\eta(z)+(A)y+(B+ns)z,y-\toto H(z)+s,z+\frac pn)$, where $A,B$ are some integers, $t\in [0,1[$, $s\in\left[0,1\right[$ and $p\in\llbracket0,n-1\rrbracket$.\medskip

If $t\not\in \Q$ then $L$ and $\eta$ are constant. Otherwise,
Let $r$ be the smallest integer such that $rt\in\N^\ast$. Again by Proposition \ref{holonomy}, $L$ is $\frac1r$-periodic. As $L(z+\frac1r)=L(z)$, ${\cal L}=1/(r\int_0^{1/r} \frac{1}{L^2(s)}\dd s)$ thus $H(z+\frac1r)=H(z)+\frac1r$ and $\eta(z+\frac1r)=\eta(z)$.

We let the reader check that $\varphi=F^{-1}\circ \varphi_0\circ F$ where  $\varphi_0$ and $F$ are defined by:
\begin{eqnarray*}
  \wdt{\varphi}_0(x,y,z)&=&(x-\frac{n\ell}{2}(z^2-z) +Ay+ (B+ns) z ,y-\ell z+s,z+t),\\
\wdt F(x,y,z)&=&(x- n(\frac{1}{2}-z)f(z),y+f(z),H(z)),   
  \end{eqnarray*}
where $f$ is the $1$-periodic function defined by:
$$f(z)=\frac{1}{A}\left(\frac{n\toto}{2}(z-H(z))^2+(B+ns) (z-H(z)-\eta(z))\right).$$

As $\varphi_0$ preserves the canonical flat affine structure of the Heisenberg group if $n\neq0$, or of $\R^3$ if $n=0$, it preserves a flat Lorentzian connection. Moreover, if $h$ is a Lorentzian metric that reads $h=2h_{1,3}(z)(\dd x-nz\dd y)\dd z+h_{2,2}(z)dy^2 +2h_{2,3}(z)\dd y\dd z+h_{3,3}(z)\dd z^2$ then it  follows from a direct computation that  $\varphi_0\in \opna{Aff}(h)\smallsetminus\opna{Is}(h)$ if and only if 
$h_{1,3}=\frac{\toto}{A}h_{2,2}$ and $h_{2,3}=\frac{h_{2,2}}{2\toto}(\frac{2\toto(B+ns)}{A} +\toto^2  +\frac{n \toto^2}{A}+ C \frac{\toto^2}{A^2}h_{2,2})$,   where $C\in\R^*$ (if $C=0$ then $\varphi_0$ is an isometry). For any $t\in[0,1]$, let $h_t$ be the metric whose connection is preserved by $\varphi_0$ that is obtained as above and satisfies $h_{t,2,2}=(1-t)h_{2,2}+t$ and $h_{t,3,3}=(1-t)h_{3,3}$. Clearly $h_0=h$ and $h_1$ is flat. Setting $h=F^{-1}{}^*g$ (it has the right form), we take  $g_t=F^*h_t$ and obtain a path of metrics having the properties announced.
\end{itemize}

All that precedes uses only that $\varphi\in\opna{Aff}(g)$ hence holds for $\varphi\in\opna{Isom}(g)$ instead of $\varphi\in\opna{Aff}(g)\smallsetminus\opna{Isom}(g)$, producing a path $(g_t)_t$ of metrics such that $\varphi\in\opna{Isom}(g_t)$ for all $t$. The only changes are in the description of the possible forms of $\varphi$, that are more constrained. In {\bf (i)}, $s\in\frac1n\Z$ instead of $\Q$, in {\bf (ii)}, $n$ must be null else the field $\partial_y+nz\partial_x$ is not Killing, and in {\bf (iv)} $\varphi_{\ell,A,B}$ must be isometric, what means that ${\cal L}$ and $k$ are in $\Lambda\Q$ and then conditions on $(\ell,A,B)$.\medskip

We are done with the first part of the theorem. 
 
Let $h$ be a metric such as in case   {\bf (iv)} of the discussion above,  with $n>0$. We recall that the leaves of $\FF$ are closed. By Lemma \ref{dynamic}, the curvature of $h$ is constant along the leaves of $\mathcal X^\perp$ and, by Proposition \ref{champs}, $\dim(\opna{Aff}(h)/\opna{Isom}(h))=1$. We choose a map $\varphi$ in $\opna{Aff}(h)$ that induces a non equicontinuous transformation of $\overline M$, an affine non isometric vector field $K$ and a parallel vector field $X$.
 Now, let us suppose that $\opna{Aff}(h)\subset\opna{Aff}(g_0)$ where $g_0$ is a flat metric, and show that $h$ is flat. This will finish the proof.

We choose coordinates $(x,y,z)$
adapted to $g_0$ i.e.\ such that $g_0$ reads $a(dx - nzdy)dz +b \dd y^2$, with $(a,b)\in \R^2$. The flow of the vector field $X$ is in $\opna{Aff}(g_0)$ and is equicontinous, therefore it is isometric and we can assume that $X=\partial_x$.
The vector field $K$ is affine and cannot preserve $g_0$ (as its flow is not equicontinuous, see \cite{Zeghib}) therefore  we can assume that 
$K=\partial_y+nz\partial_x$. As $\varphi \in\opna{Aff}(g_0)$ and $\overline\varphi$ is non equicontinuous, we can also assume that:
$$\wdt\varphi(x,y,z)=(x-\frac{n\ell}{2}(z^2-z) +Ay+B z ,y-\ell z,z),$$
for some $(\ell,A,B)\in\Z^{\ast 2}\times\Z$.

\begin{fact}
If $h$ is a metric on $\R^3/\Gamma_n$ whose connection is invariant by
the flow $\Phi_Y$ of $Y=\partial_y+nz\partial_x$, then
$\wdt h=2h_{1,3}(dx - nzdy)dz + h_{2,2}(z)dy^2 + 2h_{2,3}(z)dydz +h_{3,3}(z)dz^2$, 
with $h_{1,3}\in \R$.
\end{fact}
{\bf Proof }
The only distributions invariant by the flow being $\R\partial_x$ and
$\opna{Span}(\partial_x,\partial_y)$, they have to be parallel and
lightlike. 
Therefore $\wdt h$ reads $2h_{1,3}(x,y,z)(dx - nzdy)dz +
h_{2,2}(x,y,z)dy^2 + 2h_{2,3}(x,y,z)dydz + h_{3,3}(x,y,z)dz^2$.
Computing $\Phi_{Y}^*h$ we see that the functions $h_{1,3}$, $h_{2,2}$ and
$h_{2,3}$ are $\Phi_Y$-invariant hence depend only on $z$.
We see also that $h_{3,3}$ satisfies for any $s\in \R$:
$$h_{3,3}\circ \Phi_Y^{s}+2ns h_{1,3}=h_{3,3}+C(s)h_{1,3}^2.$$
and therefore $h_{3,3}$ is also invariant by $\Phi_Y$. It implies
immediatly that $h_{1,3}$ is in fact constant.\boxfill$\Box$
\medskip

\noindent But if $h_{1,3}$ is constant then $h_{2,2}$
and $h_{2,3}$ are also constant (it is shown at the end of point {\bf(iv)} of the discussion above), hence $h$
is flat. 
We are done.\boxfill$\Box$

Notice that during our study, we have also determined the isometries of the metrics encountered. We deduce:

\begin{cor}\label{compacite_de_isom}
Let $(M,g)$ be a compact $3$-dimensional manifold admitting a parallel lightlike field $X$, and which is not a flat torus, up to a possible 2-cover. Then $\opna{Isom}(g)$ is not compact if and only if $(M,g)$ is isometric to a metric on $\R^3/\Gamma_n$ that reads:
\begin{eqnarray*}
2\Lambda (\dd x-nz\dd y)\dd z+L^2(z)\dd y^2 + k_1 \dd y\dd z+k_2\dd z^2,
\end{eqnarray*}
with $1/\int \frac{1}{L^2(z)}\dd z\in \Lambda\Q$, ${k_1}\in {\Lambda}\Q$ and $nk_i=0$.
\end{cor}
{\bf Proof.} As $M$ is not a flat torus up to a 2-cover, it is indecomposable with a parallel vector field $X$, unique up to proportionality, see Lemma \ref{basic_pres}. It follows from Proposition \ref{prop_nonclosed} that if $X$ has non closed leaves  then the isometry group is compact. We assume now that $X$ has closed leaves. 
Let $\varphi$ be a non equicontinuous element of $\opna{Isom}(g)$. The $\varphi$ preserves $X$ and the curvature function $r$ defined in Proposition \ref{GaussBonnet}, therefore $X$ is proportional to the gradient of $r$ (otherwise $\varphi$ would be equicontinuous). It means that $r$ is constant along the leaves of $\mathcal X^\perp$.  It means that either $g$ is flat or $\mathcal X^\perp$ has closed leaves.

We suppose first that  $\mathcal X^\perp$ has closed leaves.   It  follows from the proof of Theorem \ref{theo_equivariant} that $\overline\varphi$ is also non equicontinuous. Hence, in appropriate coordinates,  the metrics reads  $2\Lambda (\dd x-nz\dd y)\dd z+L^2(z)\dd y^2 + k_1 \dd y\dd z+k_2\dd z^2$ with  $1/\int \frac{1}{L^2(z)}\dd z\in \Lambda\Q$, $nk_i=0$. 

If $n\neq 0$, a non equicontinuous isometry is obtained by composing the affine morphism $\varphi_{\ell,A,B}$ defined in \eqref{eq_philab} by an appropriate time of the flow of $Y=\partial_y+nz\partial_x$.

If $n=0$, we know, using the notation of the proof of Theorem \ref{theo_equivariant}, that  
$\wdt \varphi(x,y,z)= (x-\eta(z)+Ay+Bz,y-\toto H(z)+s,z+t)$, the function $\eta$ being determined by $g$ and $\ell$. Such a map is an isometry if and only if 
${\ell}{\mathcal L}+A\Lambda=0$ and $\ell(2k_1-A\Lambda)+2B\Lambda=0$ (see the proof of Proposition \ref{releve}).
It follows that such an isometry exists if and only if we also have $\frac{k_1}{\Lambda}\in \Q$.

 We suppose now that $\mathcal X^\perp$ has dense leaves, therefore  $g$ is flat and $n\neq 0$. 
We know that $\overline\varphi_*\oY=\oY$ as $\oY$ is linear and has an irrational slope, it implies that $\overline \varphi$ is a translation. It follows by a direct calculation that $\varphi$ preserves any frame field $(X,V,W)$ such that $g(X,V)=g(V,W)=g(W,W)=0$ and $g(V,V)=1$, it contradicts the assumption that $\varphi$ is non equicontinuous.
\boxfill$\Box$

\begin{table}[ht!]
$$\begin{array}{|c|c|}
\multicolumn{2}{l}{\text{General notation}}\\
\hline
\opna{Isom}(g)&\text{The isometry group of $(M,g)$.}\\
\hline
\opna{Aff}(g)&\text{Its affine group \ie the diffeomorphisms preserving its Levi-Civita connection.}\\
\hline
(x,y,z)&\text{Canonical coordinates of $\R^3$.}\\
\hline
\multicolumn{2}{l}{\text{Groups and diffeomorphisms used to define the manifold $M$}}\\
\hline
\tau_i&\begin{nnarray}{c}\tau_1:(x,y,z)\mapsto(x+1,y+\tau,z),\ \tau_2:(x,y,z)\mapsto(x,y+1,z),\\\tau_3:(x,y,z)\mapsto(x+r_1,y+r_2,z+1)\ \text{with }(\tau,r_1,r_2)\in(\R\smallsetminus\Q)\times\R^2.\end{nnarray}\\
\hline
\tau'_i&\begin{nnarray}{c}\tau'_1:(x,y,z)\mapsto(x+1,y,z+\tau),\ \tau'_2:(x,y,z)\mapsto(x+r_1,y+1,z+r_2),\\\tau'_3:(x,y,z)\mapsto(x,y,z+1)\ \text{with }(\tau,r_1,r_2)\in(\R\smallsetminus\Q)\times\R^2.\end{nnarray}\\
\hline
\Gamma_n&\begin{nnarray}{c}\langle\tau_x,\tau_y,\tau_{z,n}\rangle\subset\opna{Diff}(\R^3)\ \text{with }\tau_x:(x,y,z)\mapsto(x+1,y,z),\\\tau_y:(x,y,z)\mapsto(x,y+1,z)\ \text{and } \tau_{z,n}:(x,y,z)\mapsto(x+ny,y,z+1).\end{nnarray}\\
\hline
\multicolumn{2}{|c|}{\begin{nnarray}{c}\text{{\em Terminology.} If $n=0$, $\R^3/\Gamma_0=\R^3/\Z^3$ is a torus, else we call $\R^3/\Gamma_n$ a ``parabolic}\\\text{torus'', as a suspension of the parabolic automorphism $\tau_{z,n}$ of $\R^2/\Z^2$ over $\R/\Z$.}\end{nnarray}}\\
\hline
\multicolumn{2}{l}{\text{Foliations; vector and frame fields}}\\
\hline
\XX;\XX^\perp&\begin{nnarray}{c}\text{The integral foliation of the parallel lightlike vector field;}\\\text{that of its orthogonal distribution}\end{nnarray}\\
\hline
Y&\text{Projection on $\R^3/\Gamma_n$ of the $\Gamma_n$-invariant vector field $\partial_y+nz\partial_x$ of $\R^3$.}\\
\hline
\BB_{n,\theta}&\text{Projection on $\R^3/\Gamma_n$ of the $\Gamma_n$-invariant frame field $(\partial_x,Y+\theta\partial_z,\partial_z)$ of $\R^3$.}\\
\hline
\multicolumn{2}{|c|}{\begin{nnarray}{c}\text{{\em Remark.} If $n\neq0$, $\R^3/\Gamma_n$ identifies naturally with a quotient of the Heisenberg group by a}\\\text{lattice, see below Notation \ref{notation_Gamma}, then ${\cal B}_{n,\theta}$ is left invariant, equal to $(\partial_x,\partial_y+\theta\partial_z,\partial_z)$ at zero.}\end{nnarray}}\\
\hline
\multicolumn{2}{l}{\text{Numbers and functions}}\\
\hline
\begin{nnarray}{c}L,\Lambda,k,\mu\end{nnarray}&\begin{nnarray}{c}\text{Stand for real numbers or, if variables appear {\em e.g.\@} $\mu(y,z)$,}\\\text{real-valued functions; $L>0$, $\Lambda>0$.}\end{nnarray}\\
\hline
\begin{nnarray}{c}P,P'\end{nnarray}&\begin{nnarray}{c}\text{Stand for natural integers.}\end{nnarray}\\
\hline
{\cal L}&1/\int_0^11/L^2(z)\dd z\\
\hline
\multicolumn{2}{l}{\text{Diffeomorphisms of $M$ --- and two related numbers}}\\
\hline
\Phi_V^s&\text{If $V$ is any vector field, we denote by $\Phi_V^s$ its flow at time $s$.}\\
\hline
\sigma&\sigma=\Phi_Y^1\in\opna{Diff}(\R^3/\Gamma_n).\ \text{Its action on $\R^3=\wdt M$ reads: $(x,y,z)\mapsto(x+z,y,z)$.}\\
\hline
\chi&\text{If ${\cal L}=\frac pq\Lambda$ with $p,q\in\N$, $p\wedge q=1$, $\opna{Diff}(\R^3/\Gamma_n)\ni\chi=\varphi_{q,-p,0}$ defined in \eqref{eq_philab} 
p.\@ \pageref{eq_philab}.}\\
\hline
\chi',B,b&\multicolumn{1}{|l|}{\begin{nnarray}{l}\text{\ -- if $n\neq0$, $\chi'=\chi\circ\Phi_Y^s$ with $s$ making $\chi'$ isometric, see \eqref{s_et_t_pour_isometrie}.}\\\text{\ -- if $n=0$ and ${\cal L}=\frac pq\Lambda$ with $p,q\in\N$, $p\wedge q=1$, $\chi'=\chi^b\circ\sigma^{-B}$ with $(b,B)$}\\\text{\phantom{\ \ --}defined by the fact that $b$ is the smallest integer such that}\\\text{\phantom{\ \ --}$2\Lambda B=(bq)(2k+bq{\cal L})\in 2\Lambda\Z$ (\ie such that $\chi'$ is isometric, see \eqref{matricen0Fferme}).}\end{nnarray}}\\
\hline
\end{array}$$
\caption{\label{table_notations}Notation for Tables \ref{table1} and \ref{table2}, plus a remark.}
\end{table}

\begin{table}[ht!]
$$\begin{array}{|c|c|c|c|l}
\hhline{----~}
\begin{nnarray}{c}
\text{Leaves of $\XX$ and $\XX^\perp$}
\end{nnarray}&
M&
\text{Form of $g$}&
\text{${\opna{Isom}(g)}$}\\
\hhline{====~}
\begin{nnarray}{c}\text{The closure of the}\\\text{leaves of $\XX$ are}\\\text{2-tori, equal to}\\\text{the leaves of $\XX^\perp$}\end{nnarray}&
\text{$\begin{nnarray}{c}\text{\normalsize Torus}\\M\simeq\R^3/\langle\tau_1,\tau_2,\tau_3\rangle
\end{nnarray}$}
&\begin{narray}{c}\text{\small$\left(\begin{array}{ccc}0&0&\Lambda\\0&L^2(z)&0\\\Lambda&0&0\end{array}\right)$}\\\text{in }(\partial_x,\partial_y,\partial_z),L\neq{\rm Cst}.\end{narray}&
\text{$\begin{nnarray}{c}\text{The translations}\\\text{leaving $L^2(z)$}\\\text{invariant.}\\\text{\normalsize Compact.}\end{nnarray}$}
&\text{(a)}\\
\hhline{----~}
\begin{nnarray}{c}\text{The closure of the}\\\text{leaves of $\XX$ are}\\\text{2-tori, transverse to}\\\text{the leaves of $\XX^\perp$}\end{nnarray}&
\text{$\begin{nnarray}{c}\text{\normalsize Torus}\\M\simeq\R^3/\langle\tau'_1,\tau'_2,\tau'_3\rangle
\end{nnarray}$}
&\begin{narray}{c}\text{\small$\left(\begin{array}{ccc}0&0&\Lambda\\0&L^2&0\\\Lambda&0&\mu(y)\end{array}\right)$}\\\text{in }(\partial_x,\partial_y,\partial_z),\mu\neq{\rm Cst}.\end{narray}&
\begin{nnarray}{c}\text{The translations}\\\text{leaving $\mu(y)$}\\\text{invariant.}\\\text{\normalsize Compact.}\end{nnarray}
&\text{(b)}\\
\hhline{----~}
\multirow{2}{*}
  {\raisebox{-5em}{$\begin{nnarray}{c}\text{$\XX$ has closed}\\\text{leaves, $\XX^\perp$ has}\\\text{cylindrical ones}\end{nnarray}$
  }}&
\multirow{2}{*}
  {\raisebox{-5em}{$\begin{nnarray}{c}\text{Torus or}\\\text{parabolic torus}\\M\simeq\R^3/\Gamma_n\end{nnarray}$
  }}&
\begin{narray}{c}\text{\small$\left(\begin{array}{ccc}0&0&\Lambda\\0&L^2(y,z)&\nu(y,z)\\\Lambda&\nu(y,z)&\mu(y,z)\end{array}\right)$}\\\text{in }\BB_{n,\theta},\ \theta\not\in\Q\end{narray}&
\multirow{2}{*}{\raisebox{-5em}{$\begin{nnarray}{c}\text{Compact.}\end{nnarray}$}}
&\multirow{2}{*}
  {\raisebox{-5em}{(c)}}
\\
\hhline{~~-~~}
&
&\begin{nnarray}{c}\text{\small More precise form}\\\text{\small if $\theta$ is diophantine:}\\\text{\small$\left(\begin{narray}{ccc}0&0&\Lambda\\0&L^2&k\\\Lambda&k&\mu(y,z)\end{narray}\right)\text{in }\BB_{n,\theta}$.}\\{\cal J}(z):=\int_0^1\mu(y,z)\dd z\\\text{is a constant in $\left[0,2\Lambda\right[$,}\\\text{null if $n\neq0$; $k\in\left[0,\Lambda\right[$}\\\text{is null if $n\neq0$}\end{nnarray}
&\\
\hhline{----~}
\begin{nnarray}{c}\text{$\XX$ has closed}\\\text{leaves, $\XX^\perp$ has}\\\text{toric ones}\end{nnarray}&
\begin{nnarray}{c}\text{Torus or}\\\text{parabolic torus}\\M\simeq\R^3/\Gamma_n\end{nnarray}
&\begin{nnarray}{c}\text{\small$\left(\begin{narray}{ccc}0&0&\Lambda\\0&L^2(z)&k\\\Lambda&k&\mu(y,z)\end{narray}\right)$}\\\text{in }\BB_{n,0},\ \text{with the}\\\text{same conditions as}\\\text{above on ${\cal J}(z)$ and $k$.}\end{nnarray}
&
\begin{nnarray}{c}\text{Compact, except}\\\text{if ${\cal L}$ and $k$ are}\\\text{both in $\Lambda\Q$ and}\\\text{$\mu$ is constant.}\end{nnarray}
&\text{(d)}\\
\hhline{====~}
\multicolumn{4}{|c|}{\begin{nnarray}{c}\text{In cases (c)-(d), any compact group of isometries acts by translations on $\overline M=M/\XX$.}\\\text{In case it is non compact, $\opna{Isom}(g)=\langle\chi'\rangle\ltimes K$ with $K$ compact.}\end{nnarray}}\\
\hhline{----~}

\end{array}$$
\caption{\label{table1}The orientable compact indecomposable Lorentzian 3-manifolds $(M,g)$ with a parallel lightlike vector field ---~equal to $\partial_x$ in the chosen coordinates~--- and the cases where their isometry group is not compact. Flat Lorentzian tori have of course parallel lightlike vector fields but, being decomposable, {\em they are not in this table}, hence all flat metrics are here outcast in case $M$ is a torus.}
\end{table}

\section{Gathering the results}\label{tables}

Table \ref{table_notations} p.\@ \pageref{table_notations} (re)introduces the notation of this section, making it self-contained in this respect.

\begin{reminder} We call {\em indecomposable} a (pseudo\nobreakdash-)Riemannian manifold whose tangent bundle does not split into a non trivial direct sum of parallel subbundles.
\end{reminder}

\begin{table}[ht!]
$$\begin{array}{|c|c|c|l}
\hhline{---~}
\opna{Aff}(g)/\opna{Isom}(g)&\begin{nnarray}{c}\text{Possible }M\text{ and}\\\text{type of $(\XX,\XX^\perp)$}\end{nnarray}&\text{Possible }g\\
\hhline{===~}
\multirow{2}{*}{\raisebox{-3em}{\fbox{$\{[\Id]\}$}}}
&\begin{nnarray}{c}\text{Torus. The leaves of $\XX$}\\\text{are dense in 2-tori.}\end{nnarray}
&\begin{nnarray}{c}\text{Any metric of cases}\\\text{(a)-(b) of Table \ref{table1}}\end{nnarray}
&(1)\\
\hhline{~--~}
&\begin{nnarray}{c}\text{Torus or parabolic}\\\text{torus $\R^3/\Gamma_n$. The leaves}\\\text{of $\XX$ are closed, those}\\\text{of $\XX^\perp$ cylindrical.}\end{nnarray}
&\begin{nnarray}{c}\text{Any metric of case (c)}\\\text{of Table \ref{table1}, except those}\\\text{of (3) and (8) below}\end{nnarray}
&(2)\\
\hhline{---~}
\multirow{3}{*}
  {
  \raisebox{-5em}{$\begin{nnarray}{c}=\langle[\sigma]\rangle\\\simeq\fbox{$\Z$}\end{nnarray}$}
  }
&\begin{nnarray}{c}\text{Parabolic torus $\R^3/\Gamma_n$,}\\\text{$n\neq0$. The leaves}\\\text{of $\XX$ are closed, those}\\\text{of $\XX^\perp$ cylindrical.}\end{nnarray}
&\begin{nnarray}{c}\text{\small$\left(\begin{narray}{ccc}0&0&\Lambda\\0&L^2&0\\\Lambda&0&\mu(y)\end{narray}\right)$}\\\text{in }\BB_{n,\theta},\,\theta\not\in\Q,\,\mu\neq{\rm Cst}.\end{nnarray}
&(3)\\
\hhline{~--~}
&\begin{nnarray}{c}\text{Torus or parabolic}\\\text{torus $\R^3/\Gamma_n$. The leaves}\\\text{of $\XX$ are closed, those}\\\text{of $\XX^\perp$ toric.}\end{nnarray}
&\begin{nnarray}{c}\text{\small$\left(\begin{narray}{ccc}0&0&\Lambda\\0&L^2(z)&k\\\Lambda&k&\mu(y,z)\end{narray}\right)$}\text{in }\BB_{n,0},\\\text{with $k=0$ if $n\neq0$ and:}\\\text{\small$\begin{nnarray}{l}\text{$\bullet$ if $n\neq0$, $(L,\mu)$ admits no period}\\\text{$(\frac{1}{P},\frac{P'}{n}),\,P>1$, else see (5) and (9),}\\\text{$\bullet$ if $n=0$, $\mu$ depends really}\\\text{on $y$ or ${\cal L}\not\in\Lambda\Q$; else see (6)-(7).}\end{nnarray}$}\end{nnarray}
&(4)\\
\hhline{~--~}
&\multicolumn{2}{|c|}{\text{The exceptions among the metrics of case (6), see this case.}}\\
\hhline{---~}
\begin{nnarray}{c}=\langle[\sigma^v\psi^u]\rangle\simeq\fbox{$\Z$}\text{, where}\\\text{$\psi=\Phi_Y^{1/P}\circ\Phi_{\partial_z}^{P'/n}$ {\small(see $P,$}}\\\text{\small $P'$ on the right column)}\\\text{and $u$, $v\in\N$ are such that}\\\text{$un+vP=n\wedge P$}\end{nnarray}
&\begin{nnarray}{c}\text{Parabolic torus $\R^3/\Gamma_n$,}\\\text{$n\neq0$. The leaves}\\\text{of $\XX$ are closed, those}\\\text{of $\XX^\perp$ toric.}\end{nnarray}
&\begin{nnarray}{c}\text{As in (4), $\mu$ depending really on}\\\text{$y$, and $(L,\mu)$ invariant by some}\\\overline\psi:(y,z)\mapsto(y+\frac 1{P},z+\frac {P'}n).\\\text{\small$\begin{nnarray}{c}\text{We take then $P$ the largest}\\\text{integer such that such a}\\\text{$\overline\psi$ exists; $\psi$ is a lift of $\overline\psi$.}\end{nnarray}$}\end{nnarray}
&(5)\\
\hhline{---~}
\begin{nnarray}{c}=\langle[\sigma^v\chi^u]\rangle\simeq\fbox{$\Z$}\text{, where}\\\text{$u$, $v\in\N$ are such that}\\uB+vb=B\wedge b\\\hline\text{\em But: }\ =\langle[\sigma]\rangle\ \text{if}\ B=b=1.\end{nnarray}
&\begin{nnarray}{c}\text{Torus or parabolic}\\\text{torus $\R^3/\Gamma_n$. The leaves}\\\text{of $\XX$ are closed, those}\\\text{of $\XX^\perp$ toric.}\end{nnarray}&\begin{nnarray}{c}\text{As in (4), $\mu$ being}\\\text{constant, and ${\cal L}$ and}\\\text{$k$ being in $\Lambda\Q$.}\end{nnarray}
&(6)\\
\hhline{---~}
\begin{nnarray}{c}=\langle[\sigma],[\chi]\rangle\simeq\fbox{$\Z^2$}\text{, where}\\\text{$u$, $v\in\N$ are such that}\\uB+vb=B\wedge b.\end{nnarray}
&\begin{nnarray}{c}\text{Torus $\R^3/\Z^3$. The leaves}\\\text{of $\XX$ are closed, those}\\\text{of $\XX^\perp$ toric.}\end{nnarray}&\begin{nnarray}{c}\text{As in (4), $\mu$ being}\\\text{constant, ${\cal L}$ being in $\Lambda\Q$}\\\text{and $k$ in $\Lambda(\R\smallsetminus\Q)$.}\end{nnarray}
&(7)\\
\hhline{---~}
\multirow{2}{*}
  {
  $\begin{nnarray}{c}=([\Phi_{Y+\theta\partial_z}^s])_{s\in\R}\\\simeq\fbox{$\R$}\\\text{($\theta\not\in\Q$ in (8),}\\\text{$\theta=0$ in (9))}\end{nnarray}$
  }
&\begin{nnarray}{c}\text{Parabolic torus $\R^3/\Gamma_n$,}\\\text{$n\neq0$. The leaves}\\\text{of $\XX$ are closed, those}\\\text{of $\XX^\perp$ cylindrical.}\end{nnarray}
&\begin{nnarray}{c}\text{Flat: as in (3),}\\\text{with $\mu$ constant.}\end{nnarray}
&(8)\\
\hhline{~--~}
&\begin{nnarray}{c}\text{Parabolic torus $\R^3/\Gamma_n$,}\\\text{$n\neq0$. The leaves}\\\text{of $\XX$ are closed, those}\\\text{of $\XX^\perp$ toric.}\end{nnarray}
&\begin{nnarray}{c}\text{\small$\left(\begin{narray}{ccc}0&0&\Lambda\\0&L^2(z)&0\\\Lambda&0&0\end{narray}\right)$}\text{ in }\BB_{n,0}.\end{nnarray}
&(9)\\
\hhline{---~}
\end{array}$$
\caption{\label{table2}The possible forms of $\opna{Aff}(g)/\opna{Isom}(g)$ for orientable compact indecomposable Lorentzian 3-manifolds $(M,g)$ with a parallel lightlike vector field.}
\end{table}

\begin{te}\label{th_table1} Compact orientable indecomposable Lorentzian 3-manifolds $(M,g)$ with a parallel lightlike vector field are listed in Table \ref{table1}, with the cases where $\opna{Isom}(g)$ is non compact.
\end{te}
Proofs and possible details, in particular the determination of the isometry classes of these metrics for types (c) with $\theta$ diophantine, and (d), are given: in Proposition \ref{prop_nonclosed} for types (a)-(b), \ref{fibres} and \ref{Diophantine} for type (c), \ref{releve} for type (d) and Corollary \ref{compacite_de_isom} for the compactness of $\opna{Isom}(g)$.

\begin{te}\label{th_table2} If the quotient $G=\opna{Aff}(g)/\opna{Isom}(g)$ of some compact, orientable and time-orientable Lorentzian 3-manifold $(M,g)$ is non trivial then:\medskip

-- $G\simeq\Z/2\Z$,\medskip

-- or $(M,g)\simeq(\R^3,\wdt g)/\Z^3$, with $\wdt g$ flat, is a flat torus, then $G\simeq\opna{GL}_3(\Z)/(\opna{GL}_3(\Z)\cap\opna{Isom}(\wdt g))$,\medskip

-- or $(M,g)$ is indecomposable and has a parallel lightlike vector field.\medskip

\noindent Then Table \ref{table2} gives the list of the possible forms for $G$ in the third case, with its generator(s).
\end{te}
Proofs and additional details may be found in Proposition \ref{prop_nonclosed} for case (1), \ref{transfo irra} for cases (2), (3) and (8), \ref{champs} and Theorem \ref{theo_equivariant} and its proof for cases (4), (5), (6), (7) and (9), and also Corollary \ref{compacite_de_isom} to determine $\opna{Isom}(g)$ in cases (5), (6), (7).

\begin{rem}We do not treat here the case $[\opna{Aff}(g):\opna{Isom}(g)]=2$; $(M,g)$ is then decomposable and this case is very specific, see the beginning of section \ref{transfo}. 
\end{rem}

\begin{rem}If $\varphi\in\opna{Aff}(g) \smallsetminus\opna{Isom}(g)$ on some compact indecomposable Lorentzian  time-orientable  manifold   $(M,g)$ then after Lemma \ref{basic_pres} there is a constant $C_\varphi$ such that $\varphi^\ast g=g+C_\varphi X^\flat\otimes X^\flat$, where $X$ is the parallel lightlike vector field.  If $X$ is periodic, $\varphi$ has to preserve $X$ and $\varphi\mapsto C_\varphi$ defines an injective  morphism ${\cal C}:\opna{Aff}(g)/\opna{Isom}(g)\rightarrow\R$. When $\opna{dim}M=3$, the vector field is indeed  periodic and $\opna{Aff}(g)/\opna{Isom}(g)$ is always isomorphic to a subgroup of $\R$, more precisely: \medskip

-- $\opna{Im}{\cal C}=2\Lambda\Z$ when $\opna{Aff}(g)/\opna{Isom}(g)=\langle[\sigma]\rangle$,\medskip

-- $\opna{Im}{\cal C}=\frac{2\Lambda n}P(n\wedge P)\Z$ when $\opna{Aff}(g)/\opna{Isom}(g)=\langle[\sigma^v\psi^u]\rangle\simeq\Z$ \ie in case (5),\medskip

-- $\opna{Im}{\cal C}=\frac{2\Lambda n}b(B\wedge b)\Z$ when $\opna{Aff}(g)/\opna{Isom}(g)=\langle[\sigma^v\chi^u]\rangle\simeq\Z$ \ie in case (6),\medskip

-- $\opna{Im}{\cal C}=2\Lambda(\Z+\alpha\Z)$ with $\alpha\equiv \frac k\Lambda\ [\Q]$ some irrational number we do not compute, when $\opna{Aff}(g)/\opna{Isom}(g)=\langle[\sigma],[\chi]\rangle\simeq\Z^2$ \ie in case (7),\medskip

-- $\opna{Im}{\cal C}=\R$ in cases (8)-(9).\medskip

In particular, $\opna{Im}{\cal C}$ is not closed in case (7), even if $\opna{Aff}(g)$ and $\opna{Isom}(g)$ are always closed in $\opna{Diff}(M)$.
\end{rem}

\begin{rem}Tables \ref{table1} and \ref{table2} deal with {\em orientable} (it is specified) and {\em time-orientable} (the existence of the parallel lightlike field yields it) manifolds. In general, if $M$ is a compact indecomposable Lorentzian 3-manifold:\medskip

-- admitting a parallel lightlike field, it may be a quotient of order 2 of some manifold of Table \ref{table1} making it non orientable, when such quotients exist,\medskip

-- with $\afsuris>2$, it may be a quotient of order 2 or 4 of some manifold of Table \ref{table2} making it non orientable and/or time-orientable (\ie mapping $\partial_x$ on $-\partial_x$), when such quotients exist.
\end{rem}

\begin{rem}In cases (4)-(6), $\opna{Aff}(g)/\opna{Isom}(g)\simeq\Z$. However, its action is not the same. In (4) its generator $[\sigma]$ is the class of a diffeomorphism acting trivially on $\overline M=M/\XX$. In (5) it is the class of $\sigma^u\psi^v$ acting as a translation on $\overline M$ identified to $\R^2/\Z^2$ by its coordinates $(y,z)$. In (6), it is the class of  $\sigma^u\chi^v$, and the action of $\chi$ on $\overline M\simeq\R^2/\Z^2$ is homotopic to {\footnotesize$\left(\begin{array}{cc}1&p\\0&1\end{array}\right)$} if we denote ${\cal L}$ by $\frac pq\Lambda$.
\end{rem}

\section{Higher dimensional examples}\label{section_grande_dim}
In this section, we give  higher dimensional analogues  of the $3$-dimensional metrics having a big affine group. We start with a compact $m$-dimensional manifold $N$ endowed with a locally free $S^1$-action, typically $N=S^3$ with a Seifert bundle structure. Let $X$ be an infinitesimal generator of the $S^1$-action.  We choose an $X$-invariant  degenerate Riemannian metric $g_0$ of constant signature $(m-1,0)$ on $N$ such that $\ker g_0=\opna{Span}(X)$.

Let  $\wdt g$ be the Lorentzian metric on  $\R^2\times N$ given, in a decomposition $TN\times\R\partial_v\times \R \partial_u$  of the tangent bundle, by the  matrix
$$
\begin{pmatrix}
g_0 & 0  & \delta \\
0 & L^2 & 0 \\
^t \delta & 0 & 0
\end{pmatrix}.
$$
where $L$ is a function that is $X$ and $\partial_v$ invariant and $\delta$ is a matrix whose coefficients do not depend on $(u,v)$ and is such that $\wdt g(\partial_u,X)=1$.  We remark that $X$ is parallel.

Let $\Gamma$ be the group spanned by: 
$$\begin{array}{ll}
(u,v,n)\mapsto (u+1,v, \Phi_{X}^v(n))\\
(u,v,n)\mapsto (u,v+1,n)
\end{array}
$$
where $\Phi_{X}$ denotes the flow of $X$. This group preserves $\wdt g$ and acts freely, properly discontinuously
and cocompactly
on $\R^2\times N$. Moreover, the vector fields $X$ and  $Y:= \partial_v+u X$ are invariant by $\Gamma$. 

We denote by $M$ the quotient $(\R^2\times N )/\Gamma$, by $g$ the metric induced on $M$ by $\wdt g$ and by $Y$ the vector field. By a direct computation, we see that the flow of $Y$ preserves the connection induced by $g$ but does not preserve  $g$.
 
Similarly, the Levi-Civita connection of the metric $g$ on $N\times S^1$ given by:
$$\begin{pmatrix}
g_0 &   \delta \\
^t \delta &  0
\end{pmatrix}
$$
is invariant by $(u,n)\mapsto (u,\Phi_X^u(n))$ whereas $g$ is not.

Let us notice that, if $\wdt N\not\simeq \R^n$, then, according to \cite{carriere2}, the  manifolds considered above do not carry any flat Lorentzian metrics.

\hop\hop
Charles Boubel \hip
\begin{tabular}{ll}
Address: & Institut de Recherche Math\'ematique Avanc\'ee, UMR 7501 \\
& Universit\'e de Strasbourg et CNRS, \\
&7 rue Ren\'e Descartes, 67084 Strasbourg Cedex, France\\
E-mail:&{\tt charles.boubel@unistra.fr}
\end{tabular}
\hop
Pierre Mounoud\hip
\begin{tabular}{ll}
 Address: &  Institut de Math\'ematiques de Bordeaux, UMR 5251,\\
&Universit\'e de Bordeaux et CNRS\\
 &351, cours de la lib\'eration, F-33405 Talence, France\\
E-mail:&{\tt pierre.mounoud@math.u-bordeaux1.fr}
\end{tabular}

\end{document}